\documentclass[12pt]{article}
 
\usepackage{a4wide}
\usepackage{amssymb,amsmath}
\usepackage{graphicx}

 \newcommand{\num}{\mathrm{num}}
 \newcommand{\lcm}{\mathrm{lcm}}
 \newcommand{\Gal}{\mathrm{Gal}}
 \newcommand{\N}{\mathrm{norm}}
 \newcommand{\fr}{\bm}
 \newcommand{\frs}{\bms}

 \newcommand{\be}{\begin{equation}}
 \newcommand{\ee}{\end{equation}}
 \newcommand{\bea}{\begin{eqnarray}}
 \newcommand{\eea}{\end{eqnarray}}
 \newcommand{\ba}{\begin{array}}
 \newcommand{\ea}{\end{array}}
 \newcommand{\bt}{\begin{tabular}}
 \newcommand{\et}{\end{tabular}}
 \newcommand{\bd}{\begin{displaymath}}
 \newcommand{\ed}{\end{displaymath}}
 \newcommand{\PP}{{\cal P}}
 \newcommand{\bm}[1]{\boldsymbol #1}
 \newcommand{\bms}[1]{\mbox{\footnotesize \boldmath $#1$}}

 \newcommand{\NN}{\mathbb N}
 \newcommand{\QQ}{\mathbb Q}
 \newcommand{\RR}{\mathbb R}
 \newcommand{\FF}{\mathbb F}
 \newcommand{\CC}{\mathbb C}
 \newcommand{\ZZ}{\mathbb Z}

 \newcommand{\bigtimes}{\mbox{\Huge $\times$}}

\newtheorem{theorem}{Theorem}

\newtheorem{fact}{Fact}
\newtheorem{prop}{Proposition}

\begin{document}
 \bibliographystyle{unsrt}

 \parindent0pt

\begin{center}
{\Large {\bf Planar coincidences for $N$-fold symmetry}}
 \end{center}
 \vspace{3mm}
 
 \begin{center}
 {\bf {\sc Peter A.\thinspace B.\thinspace Pleasants$^{1)}$, 
     Michael Baake$^{2)}$, and Johannes Roth$^{3)}$}}
 \end{center}
\vspace{5mm}

 {\footnotesize
 \hspace*{6em} 1) Department of Mathematics, University of Queensland, \\
 \hspace*{6em} \hspace*{1em} Brisbane, QLD 4072, Australia

 \hspace*{6em} 2) Fakult\"at f\"ur Mathematik, Universit\"at Bielefeld, \\
 \hspace*{6em} \hspace*{1em} Postfach 100131, 33501 Bielefeld, Germany

 \hspace*{6em} 3) Institut f\"ur Theor.\ und Angew.\ Physik,
                  Universit\"at Stuttgart \\
 \hspace*{6em} \hspace*{1em} Pfaffenwaldring 57, 70550 Stuttgart, Germany

  $\mbox{ }$}

\vspace{3mm}
\vspace{5mm}
 
\begin{abstract}
The coincidence problem for planar patterns with $N$-fold symmetry is considered.
{}For the $N$-fold symmetric module with $N<46$, all isometries of the 
plane are classified that result in coincidences of finite index.
This is done by reformulating the problem in terms of algebraic number fields
and using prime factorization. The more complicated case $N \geq 46$ is
briefly discussed and $N=46$ is described explicitly.

The results of the coincidence problem also solve the
problem of colour lattices in two dimensions and its
natural generalization to colour modules.

\end{abstract}

\parindent15pt
\vspace{5mm}

\section{Introduction}

The concept of coincidence site lattices (CSLs) arises in the
crystallography of grain and twin boundaries \cite{Sass}.  Different
domains of a crystal do have a relationship: There is a sublattice
common to both domains across a boundary, and this is the CSL.  This
can be seen as the intersection of a perfect lattice with a rotated
copy of it where the set of points common to both forms a sublattice
of finite index, the CSL.  Up to now, CSLs have been investigated only
for special cases, for example for cubic or hexagonal crystals
\cite{grim}.  With the advent of quasicrystals infinitely many new
cases arise: quasicrystals also have grain boundaries, and one should
know the coincidence site quasilattices \cite{Warrington,WaLu}.  In a
rather different context, multiple coincidences of families of 1D
quasicrystals have been applied in constructing quasicrystals with
arbitrary symmetry (in higher dimensions) \cite{Pleasants,LuPl}.  An
application of these results was made by Rivier and Lawrence
\cite{rivier} to crystalline grain boundaries, which themselves turn
out to be quasicrystalline. This is an important example of the
relevance of a coincidence quasilattice. The experimental evidence was
provided indirectly by Sass, Tan and Balluffi in the 1970s
\cite{sass}, but beautifully by the observations of growth of
quasicrystalline grain at the grain boundary between two crystals by
Cassada, Shiflet and Poon \cite{cassada} and by Sidhom and
Portier\cite{sidhom}.  Gratias and Thalal \cite{gratias}, on the other
hand, used quasicrystal concepts in a different context to embed the
two crystal grains adjacent to a grain boundary in a higher
dimensional perfect lattice. So an extension of the CSL analysis to
more general discrete structures is desirable.

In this paper we give a unified treatment of the coincidence problem
for planar structures with general $N$-fold rotation symmetry,
extending previous \cite{Warrington,BP} and parallel \cite{Radu} work
and putting it in a more general setting.  This is what is needed for
quasicrystalline $T$-phases which are quasiperiodic in a plane and
periodically stacked in the third dimension.  Icosahedral symmetry in
3D requires different methods and will be described separately
\cite{next}.  Common to both is the necessity of an attack in two
stages: not only do we have to find the coincidence isometries (the
universal part of the problem), but also the specific modifications of
the atomic surfaces (also called windows or acceptance domains) that
are needed to describe the set of coinciding points.

In order to describe this scenario, we start with the coincidence
problem of the square lattice $\ZZ_{ }^2$.  The set of coincidence
transformations for $\ZZ_{ }^2$ forms a group, the generators of which
can be given explicitly through their connection with Gaussian
integers. Simultaneously, the so-called $\Sigma$-factor or coincidence
index can be calculated for an arbitrary CSL isometry.  Though this is
not new, the approach we use here can be generalized to quasiperiodic
planar patterns with $N$-fold symmetry.  The description of this more
general case and the tools necessary to tackle it is the main aim of
this article.

In two dimensions the classification of CSLs is the same as the
classification of colour lattices with rotational symmetry
(\cite{Schw}, Section 5.8).  In that setting the $\Sigma\/$-factor or
coincidence index is the number of colours and the different coloured
sublattices are the different cosets of the CSL in the original
lattice.  This is because, as long as the symmetry group consists of
2D rotations, all members of the symmetry group commute with the CSL
rotation, thus ensuring that the CSL is invariant under the symmetry
group.  For indecomposable groups in higher dimensions no non-trivial
orthogonal transformation commutes with all symmetries so there is no
longer this equivalence.  The only non-trivial rotation groups of 2D
lattices are $C_3$, $C_4$ and $C_6$.  The prime numbers $p\/$ for
which there exist $p\/$-colour lattices with these symmetries are
listed in \cite{Schw}, p.~76, and coincide with the sets of primes in
the denominators of the Dirichlet series given at the end of Section
\ref{examples} for the cases $n=3$ and $n=4$.  (Note that $N=2n$ for
$n$ odd and $N=n$ otherwise as will be explained later.)  {}For
non-lattices the solution of the CSL problem in 2D can be regarded as
a classification of colour modules in the plane.  An {\em $r\/$-colour
  $n\/$-module\/} is a pair of $n$-modules ($\cal M$, ${\cal M}_1$)
such that ${\cal M}_1$ has index $r\/$ in $\cal M$ and is invariant
under the symmetry group of $\cal M$, while no other coset of ${\cal
  M}_1$ has this property (see Section 3 for a definition of
$n$-modules).  The colour of a point in $\cal M$ is then determined by
its coset mod ${\cal M}_1$.  In this light, the results of Sections
\ref{CN1} and \ref{h>1} can be interpreted as finding, for each $n$,
the numbers $r\/$ for which there are $r\/$-colour $n\/$-modules and
what these $n\/$-modules are.

The paper is organized as follows. In Section 2, we review the
coincidence problem for the square lattice and formulate it in terms
of Gaussian integers.  This enables us to describe the group structure
and the coincidence indices explicitly and to introduce the concepts
needed for the generalization in Section 3.  There, the main structure
is derived with the aid of the algebraic number theory of cyclotomic
fields, followed by various explicitly worked out cases in Section 4.
They include 8-, 10- and {12-fold} symmetry, the most important cases
for quasicrystalline $T$-phases, and thus cover all cases linked to
quadratic irrationalities \cite{BJKS}.  In Section~5 we then show, in
an illustrative way, how to use the method for the eightfold symmetric
Ammann--Beenker rhombus pattern and the tenfold symmetric T\"ubingen
triangle tiling.  We give an explicit formula for the necessary
correction of the coincidence index.  In Section~6 we discuss certain
details to be dealt with for $N \geq 46$, where the variety of modules
rapidly increases, though this does not affect the generality of our
findings.  The case $N=46$ $(n=23)$ is presented in some detail.  This
is followed by some concluding remarks, while the two appendices cover
further examples (Appendix A) and proofs of technical results used in
Sections 3 and 4 (Appendix B).

\section{The square lattice: a warm-up exercise}
Let us consider the CSL problem for the square lattice $\ZZ_{ }^2$. We
focus on pure rotations first and deal with the easy extension to
reflections later.  Consider therefore a rotation (i.e., an element of
the group $SO(2) = SO(2,\RR)$) and ask for the condition that it maps
some lattice point to another one. Clearly, rotations through
multiples of $\pi/2$ do this. They form the cyclic group $C_4$ --- an
index 2 subgroup of $D_4$, the point group of $\ZZ_{ }^2$.
 
But there are more cases, as can already be seen from the growing
number of lattice points on expanding circles, summarized in the
coefficients of the theta-function of the lattice, cf.\ \cite{Conway},
\bea \label{theta1} 
  \Theta_{\ZZ^2} (x) & = & \sum_{q \in \ZZ^2}
  x_{}^{|q|^2}
  =  (\vartheta_3(x))_{ }^2    \nonumber \\
  & = & 1 + \sum_{M=1}^{\infty} r(M) x^M_{}              \\
  & = & 1 + 4 x + 4 x^2 + 4 x^4 + 8 x^5 + \ldots \nonumber 
\eea 
Here, $\vartheta_3(x) = \sum_{q \in \ZZ}\; x^{|q|^2} = 1 + 2x + 2x^4 +
2x^9 + \ldots$ is Jacobi's theta-function and $r(M)$ denotes the
number of integral solutions of the equation $a^2 + b^2 = M$, see
\cite{Hardy} for details on $r(M)$.  This number is only slowly
increasing but is unbounded, so there is an infinite number of
rotations that map one lattice point to another.

As is obvious (cf.~\cite{Lueck} and references therein), the set of
coincidence rotations (or CSL rotations) consists of all rotations
$R\/$ through angles $\varphi$ with $\sin(\varphi)=a/m$ and
$\cos(\varphi)=b/m$ rational, and hence is identical with the group
$SO(2,\QQ)$.  This requires integral solutions of the Diophantine
equation \be \label{pyth} a^2 + b^2 = m^2 \; , \ee where we need
consider only the {\em primitive} solutions, i.e., $\gcd(a,b) = 1$.
They are, of course, given by the primitive Pythagorean triples
\cite{Hardy}.  {}For a primitive solution, the set of coinciding
points forms a sublattice of $\ZZ_{ }^2$ of index $m$, whence $1/m$ is
the fraction of lattice points coinciding.  We call $m$ the {\em
  coincidence index} of $R$, denoted by $\Sigma^{}_{\ZZ^2} (R)$, or
$\Sigma(R)$ for short.  This index is often called the $\Sigma$-factor
\cite{grim,Warrington,WaLu}.

\subsubsection*{The number of CSL rotations with given index}
Without determining the rotations explicitly we can calculate their
possible indices and the number of different rotations with each index
as follows.

The number of {\em primitive} solutions of Eq.~(\ref{pyth}) can be
derived from the well-known formula (cf.\ \cite{Hardy}) 
\be
\label{numsol1} r(M) \, = \, 4(d_1(M)-d_3(M)) \; , 
\ee 
(where $d_k(M)$ counts the number of divisors of $M$ of the form
$4\ell+k$) for the total number of integer solutions of 
\be \label{2squares} a^2+b^2 \, = \, M\;.  \ee

If we write $M=2_{ }^z M_1M_3$, where $M_1$ and $M_3$ are maximal divisors of $M\/$
composed of primes congruent to 1 or 3 (mod 4), respectively, then 
Eq.~(\ref{numsol1}) can be equivalently expressed as
\be \label{numsol2}
    r(M)=\begin{cases} 4 \, d(M_1),& \text{if $M_3$ is a square},\\
                          0,& \text{otherwise}, \end{cases}
\ee
where $d(M_1)$ counts {\em all} the divisors of $M_1$.  When (as in our case) $M\/$
is a square, the first alternative in Eq.~(\ref{numsol2}) occurs.  The number of
primitive solutions, $r^*(m^2)$, of Eq.~(\ref{pyth}) can now be derived from the
``input-output" principle (cf.\ \cite{Hardy}, Thm 260) as
\be \label{io}
  r^*(m^2) \, = \, r(m^2)-\sum_p r\left(\bigl(
  \frac{m}{p}\bigr)_{ }^2\right)+
  \sum_{p,p'}r\left(\bigl(\frac{m}{pp'}\bigr)_{ }^2\right)-
  \sum_{p,p',p''}r\left(\bigl(\frac{m}{pp'p''}\bigr)_{ }^2\right) + \cdots \, ,
\ee
where $p\/$ runs through all prime factors of $m$, $pp'$ through all pairs of
distinct prime factors of $m$, and so on.  After substituting Eq.~(\ref{numsol2})
in the right hand side of Eq.~(\ref{io}) and then counting the contributions of
the factors of $m\/$ one at a time, it can be seen that
\be
  r^*(m^2)=\begin{cases}
     4 \, d^*(m),& \text{if $m\/$ has prime factors $\equiv1\;(4)$ only},\\
                            0,& \text{otherwise},\end{cases}
\ee  
where $d^*(m)$ counts the squarefree divisors of $m$.
We note that the number of CSL's (as distinct from CSL rotations) of index
$m\/$ in the square lattice is a quarter this number, since each is itself a
square lattice stabilized by the rotation group of the square (of order 4).
(Note however that not every sublattice with square symmetry is a CSL.)

So far we have:
\begin{theorem} \label{square}
The coincidence indices of the square lattice are precisely the numbers $m\/$ 
with prime factors  $\equiv 1\pmod4$ only.  The number of coincidence rotations
$\widehat{f} (m)$ with a given index $m$ is
\be
      \widehat{f}(m)=4d^*(m)
\ee
and the number of CSL's with index $m$ is
\be \label{count1}
	f(m) = d^*(m).
\ee
\end{theorem}

\subsubsection*{CSL rotations and Gaussian integers}
We have settled the question of what numbers occur as coincidence
indices of CSL rotations of $\ZZ_{ }^2$ and how many rotations there
are with each index, but there is still more to be said.

We have seen that the set of CSL rotations forms a group
($SO(2,\QQ)$, in fact). Let us introduce the notation
\be 
     SOC(\ZZ_{ }^2) := \{R \in SO(2) \; \mid \; \Sigma(R) < \infty \}
\ee
for it.  We shall investigate its structure
and derive independent generators.

The most transparent proof of Eq.~(\ref{numsol1}) (that given in \cite{Hardy})
depends on factorization in the ring of Gaussian integers.  By making 
direct use of this idea we not only find independent generators for
$SOC(\ZZ_{ }^2)$ but also have a method that readily generalizes to other
lattices and modules.

To this end, we consider the lattice $\ZZ_{ }^2$ as the ring
$\ZZ[i]$ of Gaussian integers, i.e., with $i=\sqrt{-1}$,
\be
    \ZZ[i] = \{a + i b \; \mid \; a,b \in \ZZ \}
\ee
together with the (number theoretic) norm
\be
     \N(a + i b) \, = \, (a + i b) (a - i b)
                \, = \, |a + i b|^2 \; .
\ee
The ring $\ZZ[i]$ consists of all algebraic integers in  
$\QQ(i) = \{a+ib \; \mid \; a,b \in \QQ \}$, which is both a
quadratic and a cyclotomic field. 
The coincidence rotation problem is then equivalent to finding all
numbers of norm 1 in $\QQ(i)$ because rotation through an angle means 
multiplication with the corresponding complex number on the unit circle
and a coincidence can only happen if this complex number is in $\QQ(i)$.

Any such number can uniquely be written 
(up to units) as the quotient of two Gaussian integers,
\be \label{rational1}
   e_{ }^{i \varphi} = \frac{\alpha}{\beta} = \frac{a + i b}{c + i d}
\ee
with {\em coprime} Gaussian integers $\alpha,\beta$ of identical norm,
$\N(a + i b) = \N(c + i d) = \ell$, say.
Now, we can profit from unique factorization in $\ZZ[i]$ because
any integer $\alpha \in \ZZ[i]$ divides its norm: 
\be
       \alpha \mid \N(\alpha) \; .
\ee
Let $\ell = p_1^{\nu_1} \cdots p_r^{\nu_r}$ be the (unique)
factorization of $\ell$ of (\ref{rational1}) into ``ordinary'' primes
of $\ZZ$, called {\em rational primes} from now on.  If any rational
prime $p_j$ stayed prime in $\ZZ[i]$ (i.e., did not split into two
Gaussian integers), which happens if $p_j \equiv 3$ (4), it would
appear both in the numerator and the denominator of (\ref{rational1})
which is inconsistent with coprimality.  Thus such a rational prime
cannot divide $\ell$.

A similar argument applies to the prime $2$ which, although it splits
as $2 = i (1 - i) ^2$, also has only one Gaussian prime factor up to
units.  The remaining primes are $\equiv$ 1 (mod 4) and split as $p =
\omega_p \overline{\omega}_p$ into two Gaussian integers. One of them
appears in the numerator of (\ref{rational1}), the other in the
denominator, if $p \mid \ell$.  Of course, the actual choice of
$\omega_p$ is only unique up to units and up to taking the complex
conjugate which reflects the point symmetry of the square lattice! One
convenient choice for uniqueness (which we will now take) is a
rotation angle in the interval $(0,\pi/4)$.

This, in fact, solves the above problem constructively: any CSL rotation
can be written in the form
\be \label{product1}
   e_{ }^{i \varphi} = \varepsilon \cdot \!\!\!
          \prod_{\PP \ni \, p \equiv 1 \; (4)}
         \left( \frac{\omega_p}{\overline{\omega}_p} \right)_{ }^{n_p}
\ee
where $n_p\in\ZZ$, $\varepsilon$ is a unit in $\ZZ[i]$ and $\PP$ denotes 
the set of rational primes. Since the group of units in $\ZZ[i]$ is nothing 
but $C_4$, we find
\be \label{dirprod1} 
   SOC(\ZZ_{ }^2) \simeq  C_4 \times \ZZ^{(\aleph_0)}
\ee
and the generators are $i$ (for $C_4$) and $\omega_p /
\overline{\omega}_p$ for rational primes $p \equiv 1 \; (4)$. By
$\ZZ^{(\aleph_0)}$ we mean, as usual, the infinite Abelian group that
consists of all {\em finite} integer linear combinations of the
(countably many) generators.  The coincidence index $m$ is obviously 1
for the units in $C_4$ and $p = \N(\omega_p)$ for the other generators
because this counts the number of residue classes of the CSL in
$\ZZ^2$. If the CSL rotation $R$ is factorized as in
Eq.~(\ref{product1}), we thus find 
\be \label{sigma2} \Sigma(R) =
\prod_{\PP \ni \, p \equiv 1 \; (4)} p_{ }^{|n_p|} \; .  
\ee

This solves the rotation part in principle, one can now work along the
primes $p \equiv 1 \; (4)$ to write down the generators explicitly, e.g.,
\bd 
 \frac{4+3i}{5} \; , \; \frac{12+5i}{13} \; , \; \frac{15+8i}{17} \; , \;
 \frac{21+20i}{29} \; , \; \frac{35+12i}{37} \; , \;
 \frac{40+9i}{41} \; , \; \mbox{etc ,}
\ed
where the number on the unit circle is shown in a form with denominator $p$
and rotation angle in $(0,\pi/4)$.  All other CSL rotations are obtained by
combinations, and one can regain the formula of Theorem \ref{square} for the
number of them with index $m$.
Since $d^*(m)$ is a multiplicative function
(i.e., $d^*(m_1m_2)=d^*(m_1)d^*(m_2)$ for coprime $m_1,m_2$) and $d^*(p^r)=2$ for
a prime power $p^r$ ($r \geq 1$), we obtain for $f(m)=d^*(m)$ the Dirichlet
series generating function \cite{Hardy}
\bea \label{genfun1}
\Phi(s)=\sum_{m=1}^{\infty}\frac{d^*(m)}{m^s}
&=&\prod_{p \equiv 1 \, (4)}\left(1+\frac{2}{p^s}+\frac{2}{p^{2s}}+\cdots\right) \nonumber \\
&=&\prod_{p \equiv 1 \, (4)}\frac{1+p^{-s}}{1-p^{-s}}
\eea
and the Dirichlet series generating function for $\widehat{f}(m)$ is
$4\Phi(s)$.

{}Finally, the full group of CSL isometries,
$OC(\ZZ_{ }^2)$, is the semidirect
product of the rotation part 
$SOC(\ZZ_{ }^2)$ (normal subgroup) with the
group $\ZZ_2$ generated by complex conjugation ($=$ reflection in the
$x$-axis):
\be \label{isom1}
     OC(\ZZ_{ }^2) = 
     SOC(\ZZ_{ }^2) \times_s \ZZ_2 \; .
\ee
Here conjugation of a rotation through an angle $\varphi$ by complex
conjugation results in the inverse rotation through $-\varphi$.
Let us give a brief justification of Eq.~(\ref{isom1}).
Since $O(2) = SO(2) \times_s \ZZ_2$ (semidirect product)
with the $\ZZ_2$ of Eq.~(\ref{isom1}),
{\em any} planar isometry $T$ with $\det(T)=-1$ can uniquely be written
as the product
\be \label{factor1}
     T \, = \, R(\varphi) \cdot T_x \,
\ee
of a rotation through $\varphi$ with $T_x$, the reflection in the $x$-axis.
But $T_x$ leaves $\ZZ_{ }^2$ invariant, so $T$ is a coincidence isometry
if and only if $R(\varphi)$ is a coincidence rotation.

The calculation of coincidence indices is also simple in this case.
The coincidence index for the reflection $T_x$ is 1.
{}For an arbitrary element of $OC(\ZZ_{ }^2)$, we either meet a rotation
(where we know the result already) or use the factorization (\ref{factor1}) again.
Then, the coincidence index is identical with that of its rotation part,
so Eq.~(\ref{factor1}) is all that is needed.   
This solves the coincidence problem for the square lattice completely
and we have

\begin{theorem}
The group of coincidence isometries of the square lattice $\ZZ^2$ is
\be
   OC(\ZZ^2) \, \simeq \, O(2,\QQ) \, \simeq \, 
       (C_4 \times \ZZ^{(\aleph_0)}) \, \times_s \, \ZZ_2 \; .
\ee
This group is fully characterized by Eqs.~$(\ref{product1})$, 
$(\ref{dirprod1})$ and $(\ref{isom1})$, and the coincidence index of an element
$(\ref{factor1})$ is given by Eqs.~$(\ref{sigma2})$ and $(\ref{product1})$.
\end{theorem}

\section{More generality: the unique factorization case} \label{CN1}

As briefly explained in the introduction, the corresponding programme
for a locally finite tiling ${\cal T}$ with $N$-fold symmetry (or
rather for its set of vertex sites) consists of two steps, the first
being the solution of the coincidence problem for the limit
translation module ${\cal M}({\cal T})$ of $\cal T$ (see \cite{Martin}
for details about this concept).  For the moment, we consider only
tilings with the property that the set of vertex sites of ${\cal T}$
is a subset of ${\cal M}({\cal T})$, a condition we shall come back to
in Section 5. Furthermore, we assume ${\cal M}({\cal T})$ to be what
is termed an ``$N\/$-lattice" in \cite{MRW} but which we shall call an
``$n\/$-module" (where $N=2n$ for $n$ odd and $N=n$ otherwise) in line
with the mathematical practice of reserving the word ``lattice" for
discrete subgroups.  The {\em principal $n\/$-module} (the ``standard
$N\/$-lattice" of \cite{MRW}) is the additive subgroup of $\RR_{ }^2$,
generated by the vectors of the regular $n$-star, 
\be
\label{nstar}
(\cos(2\pi k/n),\sin(2\pi k/n)) \; , \;\; k=0,\ldots,n-1 \; . 
\ee
The other modules are the non-trivial subgroups of the principal module that are
invariant under rotation about the origin through $2\pi/n$.  Modules that
differ only in scale and orientation are regarded as equivalent.

Because all modules are invariant under rotation through $\pi$ (since
if $\bm x$ is in the module then so is $-\bm x$), an $n\/$-module with
$n$ odd is invariant not only under rotation through $2\pi/n\/$ but
also through $\pi/n$.  So $n\/$-modules and $N\/$-modules are the
same. In view of this we shall assume throughout that $n\/$ is either
odd or divisible by 4, though this necessitates bearing in mind that
for odd $n\/$ an $n\/$-module has $2n\/$-fold symmetry.  The opposite
convention is used in \cite{MRW}, but the one used here is more
convenient for expressing results about cyclotomic fields that we
shall need later because it gives $n\/$ the parity of the discriminant
of the corresponding field.

The first stage of our analysis, occupying all but Section~\ref{application},
is to investigate coincidence rotations for modules and their associated
coincidence site modules, which we designate CSMs.

\subsubsection*{Symmetric modules and cyclotomic fields}

Viewed as complex numbers, the vectors (\ref{nstar}) are $\xi_{ }^k$,
where $\xi$ is a primitive $n\/$th root of 1, and the modules are
subsets of the cyclotomic field $K=\QQ(\xi)$.  The principal module is
precisely the ring of integers ${\cal O}_K$ of $K$, since it is known
that $\{1,\xi,\xi_{ }^2,\ldots,\xi_{ }^{\phi(n)-1)}\}$ is a basis for
the integers of $K$, where $\phi(n)$, the Euler totient function of
$n$, is the degree of $K\/$ over $\QQ$, cf.\ Chapter~9 of \cite{Cohn}.
The other modules are the ideals of ${\cal O}_K$ (to be defined
later), modules being equivalent precisely when they belong to the
same ideal class (defined in Section~\ref{h>1}).

In this section, at the expense of discussing only 29 modules (see
\cite{MM,MRW}), we restrict attention to values of $n\/$ for which all
$n\/$-modules are equivalent.  Because of the connection with
algebraic number theory we call this the ``class number 1" case and
use the designation ``CN1" to indicate results that are special to
this case.  (The reason behind this terminology is explained in
Section 6. Briefly, it is the case when the $n$th cyclotomic field has
class number 1.)  The class number 1 assumption simplifies the
treatment in two ways:
\begin{itemize}
  \item[1)]it is enough to solve the coincidence problem for the principal
   module
   ${\cal O}_K$ only, since all others are equivalent to it; and
  \item[2)]in the class number 1 case each integer in ${\cal O}_K$ has a
   factorization into irreducible integers that is unique apart from
   multiplying the factors by units.  (Because of the unique factorization
   these irreducible integers can safely be called {\em primes\/} in the
   class number 1 case.)
\end{itemize}

Though a convenience, the restriction to class number 1 is by no means
essential: with only minor modifications our method applies to any 2D module,
as outlined in Section \ref{h>1}.

As in the previous section, a coincidence rotation that takes $\beta$
to $\alpha$, say, ($\alpha,\beta \in {\cal O}_K$) can be represented
by the point $\gamma=\alpha/\beta$ on the unit circle.  So the CSM
problem amounts to finding the structure of the set of numbers
$\gamma$ in $K\/$ with \be \label{abs1} |\gamma|=1 \ee (a subgroup of
the multiplicative group of $K\/$).

The CSM associated with $\gamma$ is ${\cal O}_K \cap \gamma {\cal O}_K
= \num(\gamma) {\cal O}_K$, where $\num(\gamma)$, the numerator of
$\gamma$, is given by 
\be \label{den1} 
  \num(\gamma) = \gcd (\nu \in
  {\cal O}_K \; \mid \; \nu / \gamma \in {\cal O}_K ) \; , 
\ee 
and is unique up to multiplication by a unit.  In particular,
$\num(\gamma) \, \mid \, \alpha$.  The index of this module in the
original module ${\cal O}_K$ is $\N(\num(\gamma))$, the absolute norm
of $\num(\gamma)$, (\cite{Cohn} 4.4 and Cor.~2.96).  (Since units have
norm 1 this is independent of the particular numerator chosen.  All
conjugates of the field $K\/$ are complex, so norms of numbers in
$K\/$ are products of pairs of complex conjugates and hence positive.)

Eq.~(\ref{abs1}) can be reformulated as an algebraic condition with the
aid of the maximal real subfield $L\/$ of $K\/$:
\be \label{real1}
        L := \QQ(\xi+\xi_{ }^{-1})=\QQ(\cos\frac{2\pi}{n}) \; .
\ee
It is known that when $K\/$ has unique factorization $L\/$ does too
(see p.~231 of \cite{Masley2}).
As an extension of $L$, $K\/$ has degree 2 and the set of conjugates over
$L\/$ of a number $\gamma\in K\/$ is just the complex conjugate pair
$\{\gamma,\overline{\gamma}\}$.  Consequently, the {\em relative norm\/} of
$\gamma$ over $L$, $\N_{K/L}(\gamma)$, is given by
\be \label{relnorm1}
       \N_{K/L}(\gamma) = |\gamma|^2 \; .
\ee
In this notation, the {\em absolute norm} of $\gamma$ is 
$\N(\gamma) = \N_{K/\QQ}(\gamma)$ and we have the relation
\be \label{relnorm2}
       \N_{K/\QQ}(\gamma) = \N_{L/\QQ}(\N_{K/L}(\gamma)) \; .
\ee
Relative norms of integers in $K\/$ are integers in $L\/$ and norms of
units are units.  As in the previous section (where $L=\QQ$),
$\alpha\mid \N_{K/L}(\alpha) = \alpha\overline{\alpha}\/$ for every
integer $\alpha$ of $K$, so the only possible prime factors of
$\alpha$ in $K\/$ are those that divide $\N(\alpha)$.

\subsubsection*{Cyclotomic numbers on the unit circle}

When a planar module ${\cal M}$ intersects a rotated or reflected copy
of itself in a submodule of finite index, the isometry (rotation or
reflection) is again called a coincidence isometry.  The set of
coincidence isometries of ${\cal M}$ is denoted by $OC({\cal M})$. It
is again a group, with $SOC({\cal M})$ being its subgroup of
rotations.  (These concepts can be put in a much more general setting.
Some slight extensions of them are already required for the examples
in Appendix~A, for example.)

In view of Eq.~(\ref{relnorm1}) and our representation of $SOC({\cal O}_K)$ as the
elements of $K\/$ on the unit circle, we have
\be \label{norm1}
    SOC({\cal O}_K) \; \simeq \; \{ \gamma \in K 
    \; \mid \; \N_{K/L}(\gamma) = 1 \} \; .
\ee
To analyze the right hand side further we need some facts about the
arithmetic of $K\/$ and $L$.  First, the units $\varepsilon$ of $K\/$
with $|\varepsilon|=1$ are precisely the powers of $\xi$, though in
general there are also infinitely many units not on the unit circle.
(This follows, e.g., from \cite{W}, Lemma 1.6, and the last sentence
of the remark following it.)  Second, if a prime $\varrho$ of $L\/$
has two non-associated prime factors in $K\/$ (i.e., their ratio is
not a unit) then they can be taken as complex conjugates, $\omega$ and
$\overline{\omega}$.  This is because $\omega\mid\varrho$ implies
$\overline{\omega}\mid\varrho$, so, if $\omega$ and
$\overline{\omega}$ are not associates, $\omega\overline{\omega}$ is
an integer in ${\cal O}_L$ dividing $\varrho$, hence is an associate
of $\varrho$.  (Here ${\cal O}_L$ is the ring of integers of $L$, of
course.) Conversely, if $\varrho$ is divisible by just the prime
$\omega$ in $K\/$ and no other (up to units), then, as
$\overline{\omega}$ also divides $\varrho$, $\omega/\overline{\omega}$
must be a unit.  Thus a prime $\omega \in {\cal O}_K$ divides a prime
$\varrho \in {\cal O}_L$ with distinct factors if and only if
$\overline{\omega}$ is not an associate of $\omega$.  By
Eq.~(\ref{relnorm1}), $\N_{K/L}(\omega)=\N_{K/L}(\overline{\omega})$.

Now suppose that $\gamma\in K\/$ satisfies $\N_{K/L}(\gamma) = 1$ and write
$\gamma=\alpha/\beta$, where $\alpha$, $\beta$ are integers of ${\cal O}_K$ with 
no common factor.  Then
\be 
     \N_{K/L}(\alpha)=\N_{K/L}(\beta)=\nu \; \in \; {\cal O}_L
\ee
and every prime factor of $\nu$ must factorize into two non-associated
primes of $K$, one of which divides $\alpha$ only and the other $\beta$ only.
Since any such pair can be chosen to be complex conjugates, $\gamma$ can be
written as
\be \label{canon}
   \gamma=\varepsilon\prod_k \left(
   \frac{\omega_k}{\overline{\omega}_k}  \right)_{ }^{n_k},
\ee
with $\varepsilon$ a unit of $K\/$ and the $n_k$'s in $\ZZ$.  Taking
absolute values in (\ref{canon}) shows that $|\varepsilon|=1$, whence
$\varepsilon$ is a root of unity.  Different values of the $n_k$'s
give $\gamma$'s with different prime factorizations, which are
therefore not associates, and different roots of unity $\varepsilon$
give different $\gamma$'s within each set of associates.  So in this
more general situation we again have explicit presentations of
$SOC({\cal O}_K)$ and $OC({\cal O}_K)$ almost identical to those for
$SOC(\ZZ_{ }^2)$ and $OC(\ZZ_{ }^2)$ in the previous section.  These
are, for $SOC$, 
\be \label{presentation} 
SOC({\cal O}_K) \; \simeq \;
\langle\xi\rangle \!\!  \mbox{ \raisebox{-4mm}{ $\stackrel{\bigtimes}
    {\mbox{\scriptsize $\{\omega,\overline{\omega}\}\in\Omega$}}$}}
\left\langle\frac{\omega}{\overline\omega}\right\rangle \; \simeq \;
C_N\times\ZZ^{(\aleph_0)}, 
\ee 
where $\Omega$ is the set of complex conjugate pairs of non-associated
primes in $K$ (and $N=\lcm(n,2)$ as usual) and, for $OC$, 
\be
\label{OC} OC({\cal O}_K)=SOC({\cal
  O}_K)\times_s\langle\overline{\,\cdot\,}\rangle, 
\ee 
where $\overline{\,\cdot\,}$ is complex conjugation and its action on
$SOC({\cal O}_K)$ is clear.

The coincidence index of the typical rotation (\ref{canon}) of
$SOC({\cal O}_K)$ is the absolute norm of its numerator:
\be \label{index}
      \prod_k \; (\N_{K/\QQ}(\omega_k))_{ }^{|n_k|} \; .
\ee

{}For each prime pair $\{\omega,\overline\omega\}\in\Omega$, the common
value $\N(\omega)=\N(\overline\omega)$ is a rational prime power
$p^d$.  We call these prime powers the {\em basic indices} of ${\cal
  O}_K$ and the primes $p\/$ themselves the {\em complex splitting
  primes} for $K\/$ (because, in their factorization over $K$, they
contain at least one complex conjugate pair of distinct primes).  Then
(\ref{presentation}), (\ref{OC}) and (\ref{index}) show that:
\begin{prop} {\rm (CN1)} 
An integer $m\in\NN$ is a coincidence index if and only if it is a
product of basic indices.
\end{prop}

To find out what basic indices there are and count how many members of
the group $SOC({\cal O}_K)$ have given index we need to determine how
each rational prime $p\/$ factorizes in the fields $L\/$ and $K$.  It
will turn out that the basic indices are powers of distinct primes and
that whether a power of $p\/$ is a basic index (and what this power
is) depends only on the residue class of $p\/$ mod $n$.

\subsubsection*{Factorization of primes in algebraic number fields}

Before considering $K\/$ and $L\/$ specifically we describe how primes
factorize in a general algebraic number field extension $F(\alpha)\supset F$
of degree $D$.  (We shall use the standard notation $F(\alpha)/F\/$ to
denote such an extension.  A detailed account of the material in this section
can be found in Chapter~2 of \cite{Ono}).)

Let $f(x)=0$ be the minimal equation satisfied by $\alpha$ with coefficients
in $F\/$ (so the degree of $f(x)$ is $D\/$) and let ${\cal O}$ and ${\cal O}'$
be the rings of integers of $F\/$ and $F(\alpha)$.

An {\em ideal} of $\cal O$
is a subset $\fr a$ of $\cal O$ (non-empty and $\ne\{0\}$) such that
$\alpha+\beta\in{\fr a}$, $\forall\alpha,\beta\in{\fr a}$, and
$\lambda\alpha\in{\fr a}$, $\forall\lambda\in{\cal O},\alpha\in{\fr a}$.
We use the notation $(\alpha,\beta,\ldots)_{\cal O}$ to denote the smallest
ideal containing $\alpha,\beta,\ldots$ (where these are numbers in $\cal O$).
Ideals have a natural multiplication, defined by
${\fr a}{\fr b}=(\alpha\beta\mid\alpha\in{\fr a},\beta\in{\fr b})_{\cal O}$ and
$\cal O$ itself is the multiplicative identity.
There is an infinite set of {\em prime ideals} in $\cal O$ and
every ideal can be uniquely factorized into prime ideals.

Every ideal $\fr a$ in $\cal O$ extends to an ideal $({\fr a})_{{\cal O}'}$ in
${\cal O}'$, but ${\cal O}'$ also has other ideals not of this form.  For an
ideal ${\fr a}'$ of ${\cal O}'$ the {\em relative norm},
$\N_{F(\alpha)/F}({\fr a}')$, is defined as
$(\N_{F(\alpha)/F}(\alpha)\mid\alpha\in{\fr a}')_{\cal O}$ --- an ideal of
$\cal O$.  Norms are completely multiplicative
(i.e., $\N(\fr a\fr b)=\N(\fr a)\N(\fr b)$).  Let $\fr p$ be a prime ideal
in ${\cal O}$.  Then $({\fr p})_{{\cal O}'}$ factorizes into prime ideals in
${\cal O}'$ as
\be\label{ideal factors}
       ({\fr p})_{{\cal O}'}={\fr p}_1^{e_1}\ldots{\fr p}_g^{e_g} \; ,
\ee
and for each $k=1,\ldots,g$
\be
         \N_{F(\alpha)/F}^{ }\,({\fr p}_k)={\fr p}_{ }^{d_k} \; ,
\ee
where $d_k$ is the {\em residue class degree\/} of ${\fr p}_k$.
Taking norms of both sides of Eq.~(\ref{ideal factors}) shows that
\be
       d_1e_1+\cdots+d_ge_g=D,
\ee
{}For the special case of {\em normal\/} field extensions (i.e.,
extensions where $F(\alpha)$ contains not only $\alpha$ itself but
also all other roots of $f(x)=0$) we have $e_1=\cdots=e_g$ and
$d_1=\cdots=d_g$.  So in this case 
\be \label{factors}
     ({\fr p})_{{\cal O}'}=({\fr p}_1\ldots{\fr p}_g)_{ }^e \; ,
\ee
where each ${\fr p}_k$ has the same degree $d$, $d\cdot e\cdot g=D$ and
without ambiguity we can define $\deg_{F(\alpha)}^{}({\fr p}):=d$ and
$e_{F(\alpha)}^{ }({\fr p}):=e$.  Also $e_{F(\alpha)}^{ }({\fr p})>1$ only for
the finitely many primes ${\fr p}$ that divide the discriminant of the extension
$F(\alpha)/F$.  (Such primes are called {\em ramified}.)

Another special case is field extensions where $\alpha$ can be chosen so that
${\cal O}'={\cal O}[\alpha]$ (${\cal O}'$ has a
{\em simple integral basis\/} over $\cal O$).  In this case the
factorization of a prime ${\fr p}$ of $\cal O$ into prime ideals of ${\cal O}'$
mimics the factorization of $f(x)$ into irreducible factors over the finite
residue class field ${\cal O}/{\fr p}$.  So if ${\fr p}$ factorizes as in
Eq.~(\ref{ideal factors}) then
\be
    f(x)\equiv f_1(x)_{ }^{e_1}\ldots f_g(x)_{ }^{e_g}\pmod{{\fr p}}
\ee
where each $f_k$ is irreducible of degree $d_k$ and distinct $f_k$'s
correspond to distinct primes ${\fr p}_k$.  This provides a simple way of
calculating the degrees and multiplicities of the prime factors of ${\fr p}$.

The three extensions we have to deal with --- $K/\QQ$, $L/\QQ$ and
$K/L$ --- are all normal and have simple integral bases, so all the
above results apply to them.  Also relative degrees are
multiplicative: if $p\/$ is a rational prime having a prime factor
$\varrho$ in ${\cal O}_L$ which in turn has a prime factor $\omega$ in
${\cal O}_K$ then
\be \label{reldeg}
      \deg_{K/\QQ}(\omega) \, = \, \deg_{K/L}^{ }(\omega)
                                    \cdot \deg_{L/\QQ}^{ }(\varrho).
\ee
In particular, $\deg_{K/L}(\omega)$ is the same for all prime factors $\omega$
in ${\cal O}_K$ of the same rational prime $p$.

The primes of ${\cal O}_K$ in the non-associated pairs
$\{\omega,\overline\omega\}$ are precisely the unramified primes of relative
degree 1 over ${\cal O}_L$.  In view of
(\ref{reldeg}) and the normality of $K\/$ and $L$, $\Omega$ is the set of all
pairs of distinct prime factors $\{\omega,\overline\omega\}$ in ${\cal O}_K$ that
divide rational primes $p\/$ with
\be \label{eqdeg}
      \deg_L(p)=\deg_K(p) \; \; \;(=d,\mbox{ say})
\ee
and, for any such $\omega$, the absolute norm 
(cf.~(\ref{relnorm2}) above) is
\be
\N(\omega)=p_{ }^d.
\ee
We have:
\begin{prop} {\rm (CN1)} 
The complex splitting primes for $K\/$ are the rational primes that satisfy
$(\ref{eqdeg})$ and the basic indices of ${\cal O}_K$ are the powers $p^d$ of
these primes.
\end{prop}

\subsubsection*{How to calculate the CSL group and its coincidence indices}

Getting more explicit information about $SOC({\cal O}_K)$ and its coincidence
indices comes down to finding $\deg_L(p)$ and $\deg_K(p)$ for rational primes
$p$.  The following facts are sufficient to do this; we state them
here and justify them in Appendix B. 

To reiterate our notation:
$K=\QQ(\xi)$ and
$L=\QQ(\xi+\xi_{ }^{-1})$, where $\xi$ is a primitive $n\/$th root of 1, and $p\/$
is any rational prime.
\begin{fact} \label{dK}
If $p\nmid n$ then $\deg_K(p)$ is the smallest $d\in\NN$ such that $n$
divides $p_{ }^d-1$.
\end{fact}
\begin{fact} \label{dL}
If $p\nmid n$
then $\deg_L(p)$ is the smallest $d\in\NN$ such that $n$
divides at least one of $p_{ }^d+1$ or $p_{ }^d-1$.
\end{fact}
\begin{fact} \label{ramification}
(a) If $n=p^r$, for some $r$, then $p\/$ is not a complex splitting prime and
$\deg_K(p)=1$.\\
\phantom{\bf Fact 3} (b) More generally, if $n=p^rn_1$ with $p\nmid n_1$
then $p\/$ is a complex splitting prime in $K\/$ if and only if it is a complex
splitting prime in $K_1$ (the cyclotomic field of $n_1$th roots of unity).
Moreover, $\deg_K(p)=\deg_{K_1}(p)$.
\end{fact}

Although these facts alone clearly enable us to identify the complex
splitting primes and calculate their degrees and multiplicities, it is
nevertheless worth
listing some general consequences of them.\\
{\bf Remark 1}\quad These facts show that whether a prime $p\/$ is a
complex splitting prime of $n\/$ and what its degree $d\/$ is depend
only on the
residue class of $p\/$ mod $n$.\\
{\bf Remark 2}\quad Since, for a prime $\omega$ in ${\cal O}_K$
dividing $p$, $\N(\omega)=p^d$, where $d=\deg_K(p)$, Fact \ref{dK} has
the well-known consequence that $\N(\omega)\equiv1\;(\mbox{mod }n)$
for every prime $\omega$ in ${\cal O}_K$ with $\omega\nmid n$.  
In particular, every coincidence index $m\/$
with $\gcd(m,n)=1$ satisfies $m\equiv1\;(\mbox{mod }n)$.\\
{\bf Remark 3}\quad When $p\nmid n\/$ is not a
complex splitting prime, $\deg_K(p)=2\deg_L(p)$, so $\deg_K(p)$ is
even.  Facts \ref{dK}, \ref{dL} and 3(a) show that, conversely, if
$n\/$ is an odd prime power then no prime $p\/$ with $d=\deg_K(p)$
even is a complex splitting prime.  This is because if $n\mid p_{
}^d-1=(p_{ }^{d/2}-1)(p_{ }^{d/2}+1)$ but $n\nmid
p_{ }^{d/2}-1$ then, since $\gcd(p_{ }^{d/2}-1,p_{ }^{d/2}+1)=2$ (or 1
if $p=2$), $n\mid p_{ }^{d/2}+1$, so $\deg_L(p)=d/2$.

So, for $n\/$ an odd prime power, $p\/$ is a complex splitting prime if and
only if $\deg_K(p)$ is odd and $p\nmid n$, and it is unnecessary to compute
degrees over $L\/$ in this case.\\
{\bf Remark 4}\quad By Fact \ref{dK} the unramified primes with $\deg_K(p)=1$
(i.e., the primes that {\em split completely} in $K\/$) are precisely those
$\equiv1\;(\mbox{mod }n)$.  So these primes are always complex splitting primes.\\
{\bf Remark 5}\quad Facts \ref{dK} and \ref{dL} show that, for primes
$p\equiv-1\!\!\pmod n $, $\deg_K(p)=2$ and $\deg_L(p)=1$, so these primes are
never complex splitting primes.  Consequently, for every $n$, the proportion of
integers that are coincidence indices is 0.

In the next section we apply these facts and remarks to calculate coincidence
indices of specific modules.

\subsubsection*{The number of coincidences with given index}

Let $\widehat{f}(m) = N\cdot f(m)$ be the number of elements
of $SOC({\cal O}_K)$ with index $m$.    
The computational convenience of representing $\widehat{f}(m)$ this way arises
from the fact that $f(m)$ is more fundamental: it is a multiplicative function
of $m\/$ and, as for the square lattice, it counts the CSMs
with index $m$, since the rotation group of each module is the group of roots
of unity in $K\/$ and has order $N$.  In the general case, $f(m)$ cannot
be described as simply as in Eq.~(\ref{count1}), but its Dirichlet series
generating function does have a very simple expression in terms
of the $\zeta$-functions of the fields $K\/$ and $L$.  Also, for quite sizeable
individual values of the index, the number of coincidence isometries can be
calculated from (\ref{index}) and knowledge of the identity, degrees and
exponents of the complex splitting primes (or, equivalently, from the
generating function).

{}From the decomposition (\ref{presentation}) of $SOC({\cal O}_K)$ and the
function (\ref{index}) it can be seen that $f(m)$ is multiplicative,
i.e., $\gcd(m_1,m_2)=1$ implies $f(m_1m_2)=f(m_1)f(m_2)$.  This makes its
Dirichlet series
\be
\sum_{m=1}^\infty \frac{f(m)}{m^s}
\ee
a convenient tool for studying $f\/$: it can be expressed 
\cite{Hardy} as an ``Euler product"
\be
\prod_p\left(\sum_{r=1}^\infty \frac{f(p^r)}{p^{rs}}\right)\;,
\ee
with one Euler factor for each prime $p$, and the individual Euler factors are
straightforward to compute.  (The series we obtain will all be absolutely
convergent in the right half-plane Re$(s)>1$ and extendable to meromorphic
functions on the whole plane.  For using the series formally to calculate
individual values of $f\/$ these analytic properties are irrelevant, but they
play an essential r\^ole in calculating the asymptotic average value of $f$.)

Suppose the rational prime $p\/$ is divisible by the pairs
$\{\omega_1,\overline\omega_1\},\ldots,\{\omega_{g/2},\overline\omega_{g/2}\}$
of non-associated primes in $K\/$ and that each $\omega_j$ has
$\N_{K/\QQ}(\omega_j)=p^d$.  Then $f(p^k)$ is the coefficient of $p^{-ks}$ in
\be \label{Euler long}
\left(\cdots + \frac{1}{p^{2ds}}+ \frac{1}{p^{ds}}+ 1 +
\frac{1}{p^{ds}}+ \frac{1}{p^{2ds}} + \cdots\right)_{ }^{g/2},
\ee
the product of $g/2$ two-way infinite sums (one for each pair
$\{\omega_j,\overline\omega_j\}$) each having one term for each value of the
corresponding $n_k$ in (\ref{index}).  (The symmetry of the sums arises from
the fact that the index depends only on $|n_k|$, of course.)  On summing the
series this becomes
\be \label{Euler short}
\left(\frac{1+p^{-ds}}{1-p^{-ds}}\right)_{ }^{g/2}.
\ee
Since $f(m)$ is multiplicative, for a general $m\/$ it is the coefficient of
$m^{-s}$ in
\be
\prod_{{\cal C} \ni p\mid m}\left(
\frac{1+p^{-ds}}{1-p^{-ds}}\right)_{ }^{g/2},
\ee
where $\cal C$ is the set of complex splitting primes for $K$, and the values
of $d\/$ and $g\/$ are those appropriate to each individual prime $p$.  This,
in turn, is the coefficient of $m^{-s}$ in the infinite product
\be
\label{Eulerprod}
\Phi_K(s)=\prod_{p\in\cal C}\left(
\frac{1+p^{-ds}}{1-p^{-ds}}\right)_{ }^{g/2}.
\ee

To express this more simply we introduce the Dedekind $\zeta$-functions of
number fields \cite{Cohn,W}.  
The $\zeta$-function of a general algebraic number field $F\/$
is the Dirichlet series generating function for the number of ideals
$\fr a$ of $\cal O$ with $\N({\fr a})=m$, hence is given by
\be
\zeta_F(s) \, = \, \sum_{\frs a} \frac{1}{\N({\fr a})^s} \; = \;
\prod_{\frs p}\left(1- \frac{1}{\N({\fr p})^s}\right)^{-1},
\ee
where ${\fr a}$ runs through all ideals of $F\/$ and ${\fr p}$ through all prime
ideals.  When $F\/$ is normal we can collect together prime ideals ${\fr p}$
dividing the same rational prime $p\/$ to put the product on the right in the
form
\be
\prod_p\left(1- \frac{1}{p^{ds}}\right)^{-g},
\ee
where, for each rational prime $p$, $d=\deg_{F/\QQ}(p)$ and $g\/$ is the
number of prime ideals of $F\/$ dividing it.

A particular case of this is
\be
\zeta_{\QQ}(s)=\sum_{m=1}^\infty \frac{1}{m^s}=
\prod_p\left(1- \frac{1}{p^s}\right)^{-1},
\ee
which is the Riemann $\zeta$-function $\zeta(s)$ itself.

The following table compares the Euler factors of $\zeta_K(s)$ and
$\zeta_L(2s)$ for each rational prime $p$, there being three cases to
consider.  (It follows from Lemma~3 of \cite{CF} and Prop.~2.15(b) of \cite{W}
that the third case, $e_K(p)\ne e_L(p)$, occurs for at most one prime $p\/$:
the prime, if any, a power of which is equal to $n$.)
$$\bt{|c|c|c|c|c|}
\hline 
&&&&\\[-10pt]
$p$ & Field & Degree & \bt{c} Distinct \\ prime factors \et & Euler factor\\[2pt]
\hline
&&&&\\[-10pt]
\bt{l} Complex splitting\et& $\ba{c} K \\ L \ea$ & $\ba{c} d \\ d \ea$ &
          $\ba{c} g \\ g/2 \ea$ & 
          $\ba{c}(1-p^{-ds})^{-g} \\ (1-p^{-2ds})^{-g/2} \ea$ \\[2pt]
\hline
&&&&\\[-10pt]
\bt{c}Not complex splitting\\and $e_K(p)=e_L(p)$\et&
$\ba{c}K\\L\ea$&$\ba{c}d\\d/2\ea$&$\ba{c}g\\g\ea$&
$\ba{c}(1-p^{-ds})^{-g}\\(1-p^{-ds})^{-g}\ea$
\\[2pt]
\hline
&&&&\\[-10pt]
\bt{l}$e_K(p)\ne e_L(p)$\et&
$\ba{c}K\\L\ea$&$\ba{c}1\\1\ea$&$\ba{c}1\\1\ea$&
$\ba{c}(1-p^{-s})^{-1}\\(1-p^{-2s})^{-1}\ea$\\[2pt]
\hline
\et$$

On taking the quotients of the Euler factors arising from $K\/$ and $L\/$ and
comparing with Eq.~(\ref{Eulerprod}), we see that
\be
\frac{\zeta_K(s)}{\zeta_L(2s)}=\Phi_K(s)\left(1+\frac{1}{p^s}\right)^*,
\ee
the star indicating that the second factor on the right is present only if
$n\/$ is a power of a prime $p$.  So
\be \label{Dirichlet series}
CN1 \; \Rightarrow \; \Phi_K(s)=\begin{cases}
(1+p^{-s})^{-1} \zeta_K(s)/\zeta_L(2s),&
\text{if $n\/$ is a power of a prime $p\/$},\\
\zeta_K(s)/\zeta_L(2s),& \text{if not}. \end{cases}
\ee
We summarize this in
\begin{theorem}
Let $n$ be one of the $29$ numbers for which the cyclotomic field $K$ of $n$th roots
of unity has class number $1$. Then the group of coincidence rotations of an $n$-fold
symmetric module is the direct product of its finite rotation symmetry group
$C_N$ and countably many infinite cyclic groups, as in $(\ref{presentation})$, and the full
group of coincidence isometries is the extension of this by a reflection symmetry.
The coincidence index of such an isometry is given by $(\ref{index})$ and $(\ref{canon})$. 
The Dirichlet series generating function for
\begin{eqnarray}
    f(m) \,&=&\, \{\mbox{number of CSMs of index $m$}\}\\
\,&=&\,\frac{1}{N}\times\{\mbox{number of coincidence rotations of index $m$}\}
\end{eqnarray}
is given by $(\ref{Dirichlet series})$.
\end{theorem}

A principal use of a Dirichlet series is to find asymptotic formul{\ae} for sum
functions of its coefficients by means of residue calculus.  In the present
instance this technique shows, for example, that
\begin{eqnarray}
\mbox{Number of CSMs of index $<X$}&=&\sum_{m<X}f(m)\nonumber\\
&\sim&X\cdot\{\mbox{residue of $\Phi_K(s)$ at $s=1$}\}.
\end{eqnarray}
In view of Eq.~(\ref{Dirichlet series}), this residue can be computed from
known formul{\ae} for the residues of $\zeta$-functions at 1 and values of
$\zeta$-functions at 2.  The value of the residue can be regarded as the
``average number of CSMs" with a given arbitrarily chosen
positive integer as index.

\section{Examples: $N \, =$ 6, 4, 10, 14, 8, and 12} \label{examples}

After the general derivation of the previous section, let us present some
examples explicitly.  We select those relevant to known crystals and
quasicrystals.  For each example we list:
\begin{itemize}
\item[(a)]The fields $K\/$ and $L\/$ and the degree $[K:\QQ]$ of
          $K\/$ over $\QQ$.

\item[(b)]A table giving, for each residue class mod $n\/$ containing primes
          $p$, $\deg_K(p)$ (and, if necessary, $\deg_L(p)$ too).  In the bottom
          line of the table (where $\deg_K(p)$ is given) the degrees of complex
          splitting primes are underlined.  With each table is a comment
          describing which facts and remarks from the previous section were used
          to compute it.
\item[(c)]A list of the types of basic indices, using the notation that
          $p_{(a)}^b$ represents the $b$\/th powers of all primes congruent to
          $a\/$ mod $n$.
\item[(d)]The Dirichlet series generating function of $f(m)$, given as a
          ratio of $\zeta$-functions, as an Euler product and expanded
          explicitly as far as the 12th nonzero term.  (The same notation as in
          (c) is used for the primes in the Euler product.)
\item[(e)]An explicit formula for $f(m)$ in the style of Eq.~(\ref{count1}).  In
          these formul{\ae} $e_p$ denotes the largest exponent $e\/$ for which
          $p^e\mid m$.
\item[(f)]The average value of $f(m)$, as defined above.
\end{itemize}

The smallest coincidence indices can be read off as the denominators (with
$s=1$) of the Dirichlet series, with the corresponding values of $f(m)$ as the
numerators.  All values of $f(m)$ for $m>1$ are even, reflecting the
geometrical fact that the reverse of a coincidence rotation is also a
coincidence rotation.

\subsection*{$n=3$, the triangular (or hexagonal) lattice}

\hspace{15pt} $K=\QQ(\sqrt{-3})$, $L=\QQ$, $[K:\QQ]=2$.
\bd
\ba{|r|ccc|}
\hline  
&&&\\[-10pt]
p \!\! \pmod{3} &1&2&3\\[2pt]
\hline
&&&\\[-10pt]
\deg_K(p)&\underline1&2&1\\[2pt]
\hline
\ea
\qquad
\mbox{\bt{l}Computed using\\Facts \ref{dK} and 3(a) and Remark 3.\et}
\ed

{\bf Basic indices}: $p_{(1)}^{ }$ 

{\bf Dirichlet series}:
\begin{eqnarray*}
\lefteqn{\left(1+\frac{1}{3^s}\right)^{-1}\frac{\zeta_K(s)}{\zeta(2s)} \; = \;
    \prod \frac{1 + p_{(1)}^{-s}}{1 - p_{(1)}^{-s}}}\\
& & \mbox{\small $
  =1+\frac{2}{7^s}+\frac{2}{13^s}+\frac{2}{19^s}+\frac{2}{31^s}+
    \frac{2}{37^s}+\frac{2}{43^s}+\frac{2}{49^s}+\frac{2}{61^s}+
    \frac{2}{67^s}+\frac{2}{73^s}+\frac{2}{79^s}+\cdots $ }
\end{eqnarray*}

{\bf Number of CSLs with index $m\/$:}
$$f(m)=\begin{cases}
       \displaystyle\prod_{p\mid m}2,&
       \text{if $m$ is a product of basic indices},\\
        0,& \text{otherwise}.\end{cases}$$

{\bf Average number of CSLs:}
\[
   \frac{\sqrt3}{2\pi}\simeq0.276
\]

\subsection*{$n=4$, the square lattice}

\hspace{15pt} $K=\QQ(i)$, $L=\QQ$, $[K:\QQ]=2$.
\bd
\ba{|r|ccc|}
\hline
&&&\\[-10pt]
p \!\! \pmod{4} &1&2&3\\[2pt]
\hline
&&&\\[-10pt]
\deg_L(p)&1&1&1\\
\deg_K(p)&\underline1&1&2\\[2pt]
\hline
\ea
\qquad
\mbox{\bt{l}Computed from\\Facts \ref{dK}, \ref{dL} and 3(a).\et}
\ed

{\bf Basic indices}: $p_{(1)}^{ }$

{\bf Dirichlet series}:
\begin{eqnarray*}
\lefteqn{\left(1+\frac{1}{2^s}\right)^{-1} \frac{\zeta_K(s)}{\zeta(2s)}=
     \prod \frac{1 + p_{(1)}^{-s}}{1 - p_{(1)}^{-s}}}\\
& & \mbox{\small $
   =1+\frac{2}{5^s}+\frac{2}{13^s}+\frac{2}{17^s}+\frac{2}{25^s}+
      \frac{2}{29^s}+\frac{2}{37^s}+\frac{2}{41^s}+\frac{2}{53^s}+
      \frac{2}{61^s}+\frac{4}{65^s}+\frac{2}{73^s}+\cdots $}
\end{eqnarray*}

{\bf Number of CSLs with index $m\/$:}
\[
    f(m)=\begin{cases}
       \displaystyle\prod_{p\mid m}2,&
        \text{if $m$ is a product of basic indices}, \\
        0,& \text{otherwise}.\end{cases}
\]

{\bf Average number of CSLs:}
\[  \frac{1}{\pi}\simeq 0.318  \]

\subsection*{$n=5$, the $10$-fold module}

\hspace{15pt} $K=\QQ(e^{2\pi i/5})$, $L=\QQ(\sqrt5)$, $[K:\QQ]=4$.
\bd
\ba{|r|ccccc|}
\hline
&&&&&\\[-10pt]
p \!\! \pmod{5} &1&2&3&4&5\\[2pt]
\hline
&&&&&\\[-10pt]
\deg_K(p)&\underline1&4&4&2&1\\[2pt]
\hline
\ea
\qquad
\mbox{\bt{l}Computed using\\Facts \ref{dK} and 3(a) and Remark 3.\et}
\ed

{\bf Basic indices}: $p_{(1)}^{ }$

{\bf Dirichlet series}:
\begin{eqnarray*}
\lefteqn{\left(1+ \frac{1}{5^s}\right)^{-1} \frac{\zeta_K(s)}{\zeta_L(2s)}=
    \prod \left(\frac{1 + p_{(1)}^{-s}}{1 - p_{(1)}^{-s}}\right)^2}\\
& & \mbox{\small $
   =1+\frac{4}{11^s}+\frac{4}{31^s}+\frac{4}{41^s}+\frac{4}{61^s}+
      \frac{4}{71^s}+\frac{4}{101^s}+\frac{8}{121^s}+\frac{4}{131^s}+
      \frac{4}{151^s}+\frac{4}{181^s}+\frac{4}{191^s}+\cdots $}
\end{eqnarray*}

{\bf Number of CSMs with index $m\/$:}
\[ 
    f(m)=\begin{cases}
       \displaystyle\prod_{p\mid m}4e_p,&
         \text{if $m$ is a product of basic indices}, \\
        0,& \text{otherwise}.\end{cases}
\]

{\bf Average number of CSMs:}
\[
  \frac{5\log\tau}{\pi^2}\simeq 0.244
\]

\subsection*{$n=7$, the $14$-fold module}

\hspace{15pt} $K=\QQ(e^{2\pi i/7})$, $L=\QQ(\cos(2\pi/7))$, $[K:\QQ]=6$.
\bd
\ba{|r|ccccccc|}
\hline
&&&&&&&\\[-10pt]
p \!\! \pmod{7} &1&2&3&4&5&6&7\\[2pt]
\hline
&&&&&&&\\[-10pt]
\deg_K(p)&\underline1&\underline3&6&\underline3&6&2&1\\[2pt]
\hline
\ea
\qquad
\mbox{\bt{l}Computed using\\Facts \ref{dK} and 3(a) and Remark 3.\et}
\ed

{\bf Basic indices}: $p_{(1)}^{ }$, $p_{(2)}^3$, $p_{(4)}^3$

{\bf Dirichlet series}:
\begin{eqnarray*}
\lefteqn{\left(1+ \frac{1}{7^s}\right)^{-1} \frac{\zeta_K(s)}{\zeta_L(2s)}=
    \prod \left(\frac{1 + p_{(1)}^{-s}}{1 - p_{(1)}^{-s}}\right)^3
          \frac{(1 + p_{(2)}^{-3s})(1 + p_{(4)}^{-3s})}	
               {(1 - p_{(2)}^{-3s})(1 - p_{(4)}^{-3s})}} \\
& & \mbox{\small $
    =1+\frac{2}{8^s}+\frac{6}{29^s}+\frac{6}{43^s}+\frac{2}{64^s}+
       \frac{6}{71^s}+\frac{6}{113^s}+\frac{6}{127^s}+\frac{6}{197^s}+
       \frac{6}{211^s}+\frac{12}{232^s}+\frac{6}{239^s}+\cdots $}
\end{eqnarray*}

{\bf Number of CSMs with index $m\/$:}
\[
   f(m)=\begin{cases} \displaystyle
          \prod_{\substack{p\mid m \\ p\equiv 1\; (7)}}(4e_p^2+2)
	     \prod_{\substack{p\mid m\\ p\not\equiv 1\; (7)}}2,
	     & \text{if $m$ is a product of basic indices}, \\
	     \;\;\;0,& \text{otherwise}.\end{cases}
\]

{\bf Average number of CSMs:}
\[
    \frac{21\sqrt7R}{16\pi^3}\simeq 0.235,
\]
where $R\/$ (the regulator of $K\/$) is given by
\[
   \frac{R}{4}=\log^2\bigl(2\cos\frac{2\pi}{7}\bigr)-
   \log\bigl(2\cos\frac{\pi}{7}\bigr)\log
   \bigl(2\cos\frac{3\pi}{7}\bigr)\simeq 0.525.
\]

\subsection*{$n=8$, the $8$-fold module}

\hspace{15pt} $K=\QQ(e^{\pi i/4})$, $L=\QQ(\sqrt2)$, $[K:\QQ]=4$.
\bd
\ba{|r|ccccc|}
\hline
&&&&&\\[-10pt]
p \!\! \pmod{8} &1&2&3&5&7\\[2pt]
\hline
&&&&&\\[-10pt]
\deg_L(p)&1&1&2&2&1\\
\deg_K(p)&\underline1&1&\underline2&\underline2&2\\[2pt]
\hline
\ea
\qquad
\mbox{\bt{l}Computed from\\Facts \ref{dK}, \ref{dL} and 3(a).\et}
\ed

{\bf Basic indices}: $p_{(1)}^{ }$, $p_{(3)}^2$, $p_{(5)}^2$

{\bf Dirichlet series}:
\begin{eqnarray*}
\lefteqn{\left(1+ \frac{1}{2^s}\right)^{-1} \frac{\zeta_K(s)}{\zeta_L(2s)}=
    \prod \left(\frac{1 + p_{(1)}^{-s}}{1 - p_{(1)}^{-s}}\right)^2
        \frac{(1 + p_{(3)}^{-2s})(1 + p_{(5)}^{-2s})}
	{(1 - p_{(3)}^{-2s})(1 - p_{(5)}^{-2s})}}\\
& & \mbox{\small $
    =1+\frac{2}{9^s}+\frac{4}{17^s}+\frac{2}{25^s}+\frac{4}{41^s}+
       \frac{4}{73^s}+\frac{2}{81^s}+\frac{4}{89^s}+\frac{4}{97^s}+
       \frac{4}{113^s}+\frac{2}{121^s}+\frac{4}{137^s} +\cdots $}
\end{eqnarray*}

{\bf Number of CSMs with index $m\/$:}
\[
   f(m)= \begin{cases}\displaystyle
  \prod_{\substack{p\mid m\\ p\equiv1\;(8)}}
        4e_p\prod_{\substack{p\mid m\\ p\not\equiv1\;(8)}}2
				,&
   \text{if $m$ is a product of basic indices}, \\
       0,&  \text{otherwise}.\end{cases}
\]

{\bf Average number of CSMs:}
\[
   \frac{2\sqrt2\log(1+\sqrt2)}{\pi^2}\simeq 0.253
\]

\subsection*{$n=12$, the $12$-fold module}

\hspace{15pt} $K=\QQ(e^{\pi i/6})$, $L=\QQ(\sqrt3)$, $[K:\QQ]=4$.
\bd
\ba{|r|cccccc|}
\hline
&&&&&&\\[-10pt]
p \!\! \pmod{12}&1&2&3&5&7&11\\[2pt]
\hline
&&&&&&\\[-10pt]
\deg_L(p)&1&1&1&2&2&1\\
\deg_K(p)&\underline1&2&2&\underline2&\underline2&2\\[2pt]
\hline 
\ea
\qquad
\mbox{\bt{l}Computed from Facts 1 and 2\\
	and the cases $n=3$ and $n=4$\\
	using Fact 3(b).\et}
\ed

{\bf Basic indices}: $p_{(1)}^{ }$, $p_{(5)}^2$, $p_{(7)}^2$

{\bf Dirichlet series}:
\begin{eqnarray*}
\lefteqn{
    \frac{\zeta_K(s)}{\zeta_L(2s)} =
    \prod \left( \frac { 1 + p_{(1)}^{-s} } { 1 - p_{(1)}^{-s} } \right)^2
         \frac{(1 + p_{(5)}^{-2s})(1 + p_{(7)}^{-2s})}
	      {(1 - p_{(5)}^{-2s})(1 - p_{(7)}^{-2s})} }\\
& & \mbox{\small $
   =1+\frac{4}{13^s} +\frac{2}{25^s} +\frac{4}{37^s} +\frac{2}{49^s}+
      \frac{4}{61^s} +\frac{4}{73^s }+\frac{4}{97^s} +\frac{4}{109^s}+
      \frac{4}{157^s}+\frac{8}{169^s}+\frac{4}{181^s}+\cdots $ }
\end{eqnarray*}

{\bf Number of CSMs with index $m\/$:}
\[
   f(m)= \begin{cases}\displaystyle
 \prod_{\substack{p\mid m\\ p\equiv 1\;(12)}}
     4e_p\prod_{\substack{p\mid m\\ p\not\equiv 1\;(12)}}2
				,&
 \text{if $m$ is a product of basic indices}, \\
       0,& \text{otherwise}.\end{cases}
\]

{\bf Average number of CSMs:}
\[
   \frac{\sqrt3\log(2+\sqrt3)}{\pi^2}\simeq 0.231
\]

\section{Application to 2D quasicrystals}\label{application}

The reader might like to see at least one or two examples where we apply
the above results to planar quasicrystals.
{}For simplicity, we consider the eightfold symmetric Ammann--Beenker tiling
and the decagonal T\"ubingen triangle tiling \cite{BKSZ} here, 
while the slightly more
complicated rhombic Penrose tiling is discussed in Appendix A.

\subsubsection*{The Ammann--Beenker tiling}

Consider the eightfold symmetric Ammann--Beenker tiling
of Fig.~1 and, in particular, the coincidence problem of its vertex points
for rotations around the symmetry centre. The underlying module is the
standard eightfold module of rank 4, usually obtained as projection of the
hypercubic lattice $\ZZ^4$ to a suitably chosen 2D plane. This plane,
and its perpendicular complement, are eigenspaces of an eightfold rotation.

\begin{figure}
\centerline{\includegraphics{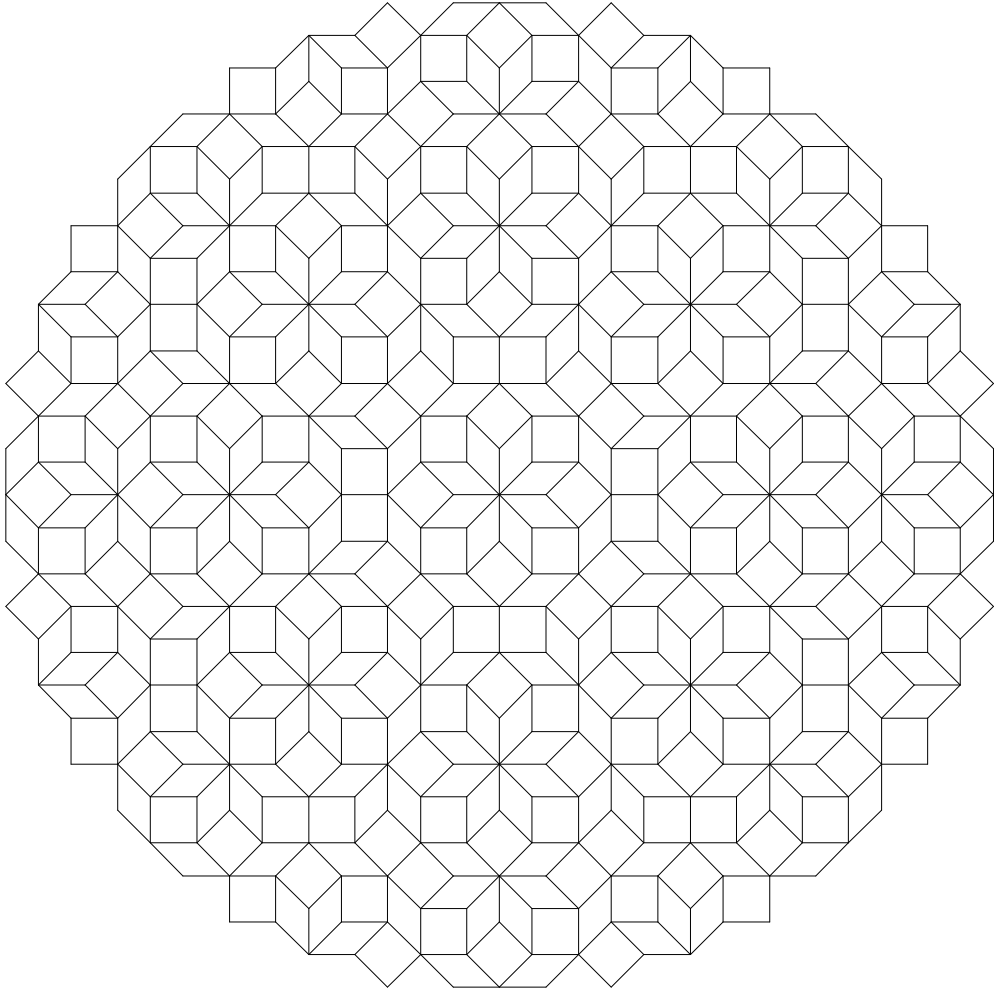}}
\caption{Central patch of the exactly eightfold symmetric
         Ammann--Beenker tiling.}
\label{fig:ab}
\end{figure}

The set of vertex sites of this tiling is just the subset of module points
whose corresponding points in $\ZZ^4$ 
perpendicularly project into a certain regular octagonal window. 
It is clear then that a coincidence of vertex sites
implies one in the module, but also the converse is true due to the way
the tiling sites are distributed over the module.

\begin{figure}
\centerline{\includegraphics{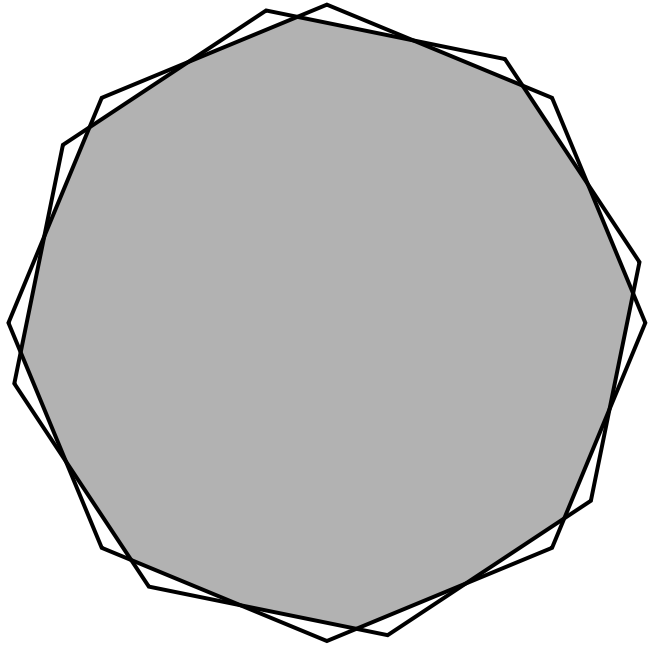}}
\caption{Intersection of two acceptance domains that are
rotated against each other.}
\label{fig:win}
\end{figure}

A coincidence rotation can be lifted to 4-space whence it also affects
the window. In fact, a coincidence point must have perpendicular projections
both in the original and in the rotated window! But this results in a slight
modification of the fraction of coinciding points which has to be corrected by an
acceptance factor $A$. This is nothing but the area ratio of the
intersection of the rotated windows with the original window, see Fig.~2.
{}For a coincidence rotation through $\varphi$, it turns out to be
\be \label{accept1}
  A = 1 - ( 1 - \frac{1}{\sqrt{2}}) \sin(\widehat{\psi}) 
           \sin(\frac{\pi}{4} - \widehat{\psi}) \; ,
\ee
where $\widehat{\psi} \in [0,\pi/4)$ via
\be
   \widehat{\psi} = \psi - \left[\frac{4\psi}{\pi}\right] 
                    \cdot\frac{\pi}{ 4}
\ee
and $\psi$, the rotation angle in perpendicular or internal space, is related to
the angle $\varphi = 2\arctan(a + b \sqrt{2})$ through an algebraic conjugation:
\be 
      \psi = 2\arctan(a - b \sqrt{2}) \; .
\ee

The acceptance factor (\ref{accept1}) is 1 for symmetry rotations 
and smaller otherwise, the minimum value
being $A_{\min} \simeq 0.957$ at $\pi/8$. 
The set of coinciding points almost looks like an Ammann--Beenker pattern
again, but some points are missing: the quantity $1-A$ is the frequency
of such failures which were observed in \cite{Warrington}.
With a more complicated window, star-shaped say,
the acceptance factor would also become more complicated: with some choices of
window it can even be zero for certain angles.
But we will not go into further details here.

\subsubsection*{The T\"ubingen triangle tiling}

Let us now consider the coincidence problem for the vertices of
the decagonal triangular tiling of Fig.~3. All vertex sites
belong to the standard tenfold module
which can be obtained by projection of the root lattice $A_4$
to a suitably chosen plane \cite{BKSZ}. For simplicity, we consider
the cartwheel tiling (which is singular) because it has full $D_{10}$ symmetry
in the sense that
a $D_{10}$ operation produces mismatches of density zero in the plane
(along worms). We thus have coincidence fraction $1$ in this case.
Also, all other coincidences of the tenfold module are realized.
As in the previous example, one has to correct the coincidence fraction,
this time by rotating a decagon (the window of the vertex sites)
and intersecting it with the original one. Let us give the
correction formula in slightly more generality.  If the window were a regular
$n\/$-gon, the analogue of Eq.~(\ref{accept1}) would read
\be \label{accept2}
  A = 1 - \left( \frac{\sin(\alpha/2)}{\sin(\alpha)} \right)^2
               \sin(\widehat{\psi}) \sin(\alpha - \widehat{\psi}) \; ,
\ee
where $\alpha = 2 \pi  /n$, $\widehat{\psi} = \psi - [ \frac{n \psi}{2 \pi}] \cdot
\frac{2 \pi}{n}$, and $\psi$ is related to $\varphi$ via an algebraic
conjugation.

\begin{figure}
\centerline{\includegraphics{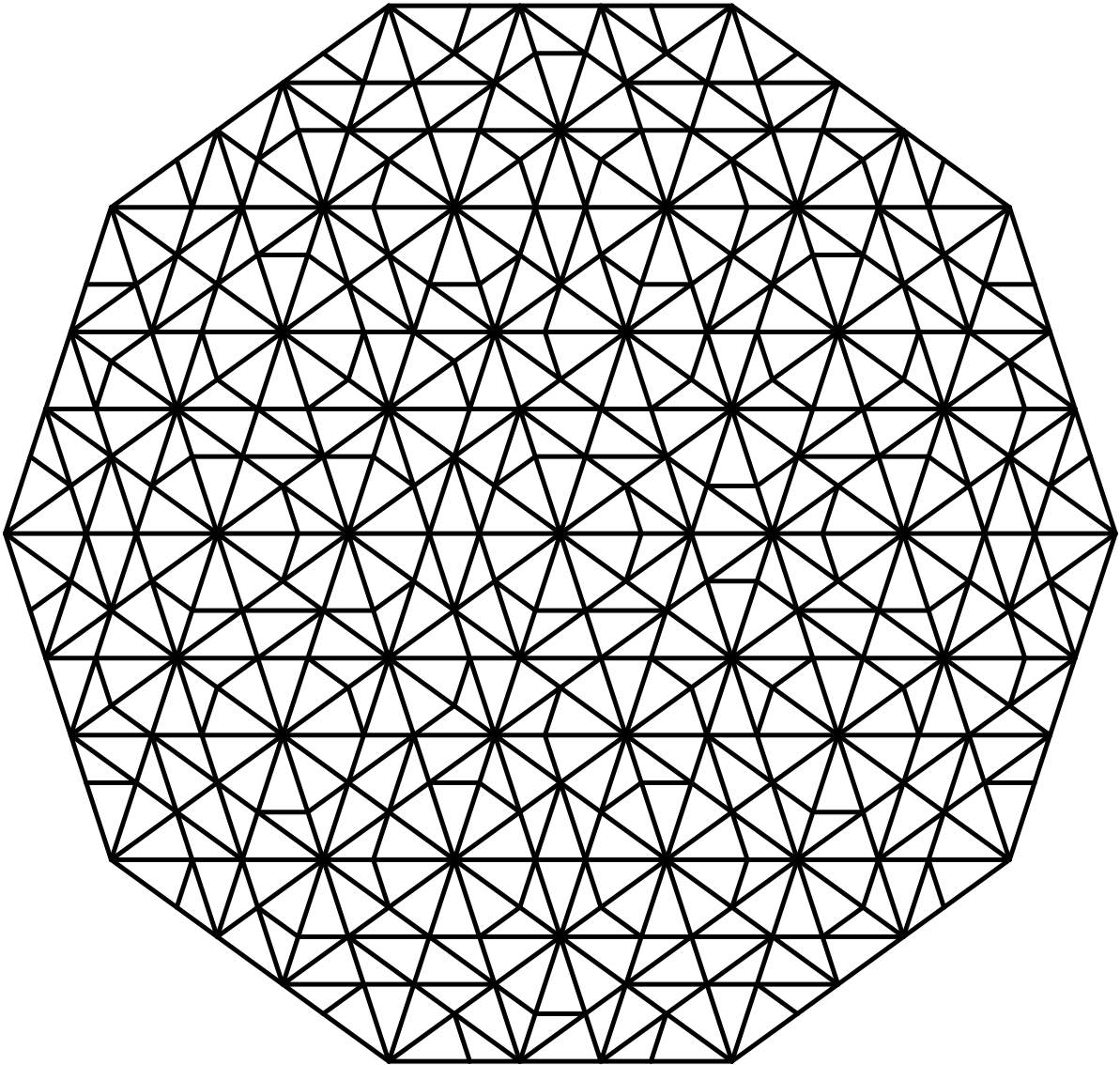}}
\caption{Central patch of the cartwheel version of the
tenfold symmetric T\"ubingen triangle tiling.}
\label{fig:ttt}
\end{figure}

In the present instance, this relation is that $\tan(\varphi/2)$ can be expressed
in the form
\be
\tan\frac{\varphi}{2}=(a+b\tau)\sin\frac{2\pi}{5}\quad(a,b\in\QQ)
\ee
and then
\be
\tan\frac{\psi}{2}=(a+b \tau')\sin\frac{4\pi}{5},
\ee
where $\tau' =-1/\tau$ is the conjugate of $\tau$ in $\QQ(\tau)$.

\subsubsection*{Twelvefold symmetric tilings}

As well as eight- and tenfold symmetries, twelvefold symmetry is of practical
interest.  Here the calculation of $\psi$ is very similar to the eightfold
case: given the angle
\be
\phi=2\arctan(a+b\sqrt3)
\ee
in tiling space one obtains the angle
\be
\psi=2\arctan(a-b\sqrt3)
\ee
in internal space, which can be used with Eq.~(\ref{accept2}).

\section{Beyond unique factorization}  \label{h>1}

In Section~3, we restricted ourselves to the ``class number 1" case, where
there is essentially only one $n\/$-module.  We now show how our method can be
adapted to other cases, too. The smallest value of $n\/$ to which Section~3
does not apply is 23 ($N=46$), mentioned in \cite{MRW}. Here, the
cyclotomic field has class number 3, so there are 3 distinct modules with
46-fold symmetry. (The number of modules increases rapidly with $n\/$
\cite{Masley1,MRW}.)

\subsubsection{Ideals and ideal classes}

Let $F\/$ be algebraic number field with ring of integers ${\cal O}$.
The set of ideals of ${\cal O}$ can be extended to form a group by admitting 
{\em fractional ideals\/} of the form
\be
{\fr a}{\fr b}^{-1}=\{\gamma\mid\gamma\beta\in{\fr a}\;\;\forall\beta\in{\fr b}\} \, ,
\ee
where ${\fr a}$ and ${\fr b}$ are ideals as defined in Section~3.
(A fractional ideal need not be a subset of ${\cal O}$.)  The identity element
of the group of fractional ideals is ${\cal O}$.  A {\em principal ideal\/} is a
fractional ideal of the form
\be
(\gamma)_{\cal O}=\{\gamma\alpha\mid\alpha\in{\cal O}\},
\ee
generated by the single number $\gamma$.  When unique factorization into
irreducible integers fails in $F\/$ then some ideals must necessarily be
non-principal.  Two fractional ideals ${\fr a}$ and ${\fr b}$ are
{\em equivalent\/} if ${\fr b}=\gamma{\fr a}$ for some $\gamma\in F$.  The
equivalence classes, called {\em ideal classes}, form a quotient group of the
group of fractional ideals called the {\em ideal class group}, $H=H(F)$, which
turns out to be finite. Its order is called the {\em class number},
$h(F)$.  The identity element of $H\/$ is the class of principal ideals.

The ideal classes inherit complex conjugation from $F\/$: each ideal class
$C\in H\/$ has a complex conjugate class $\overline C$.  Complex conjugation is
an automorphism of $H\/$ of order 2.

\subsubsection{Ideals as modules}

Our definition of $n\/$-modules makes them ideals in the ring of integers of
the $n\/$th cyclotomic field (and with any broader definition an $n\/$-module
would certainly be equivalent to one of these).  Multiplication by a complex
number $\gamma$ is equivalent to a combined rotation and scale change in the
plane, so equivalent ideals certainly correspond to equivalent modules. 
Conversely, equivalent modules can be transformed into each other by
multiplication by a complex number $\gamma$, and if both modules are subsets of
an algebraic number field $K\/$ then $\gamma$ is in $K\/$ and the corresponding
ideals are equivalent.

So the set of $n\/$-modules up to equivalence corresponds to the class group of
the $n\/$th cyclotomic field.

\subsubsection{Coincidence rotations in the general case}

With class number $>1$, $n\/$-modules are no longer all equivalent.  So, for
comprehensiveness, we need to consider not just $OC({\cal O}_K)$ but also
$OC({\fr c})$ for an arbitrary ideal $\fr c$ of ${\cal O}_K$.

There are two problems to be overcome in extending our method to the general
$n\/$-module:
\begin{itemize}
\item[(1)]How to classify which of the products on the right of
Eq.~(\ref{canon}) give rise to numbers $\gamma$ with $|\gamma|=1$ (when some
$\omega_k$'s are non-principal ideals) and
\item[(2)]how to choose a representative of the reflection coset of $OC$, for
modules not invariant under complex conjugation, and how to calculate
coincidence indices of reflections from it.
\end{itemize}
In this subsection we address the first of these.

For the fractional ideal ${\fr a}\overline{\fr a}^{-1}$ to give rise to a number
$\gamma\in K\/$ with $|\gamma|=1$ two conditions are necessary (and the
conjunction of these conditions is also sufficient).  They are
\begin{itemize}
\item[(A)]${\fr a}\overline{\fr a}^{-1}$ is principal, and
\item[(B)]for every $\delta$ such that
$(\delta)_{\cal O}={\fr a}\overline{\fr a}^{-1}$,
$\delta\overline{\delta}=\varepsilon\overline{\varepsilon}$ for some
unit
$\epsilon$ of $K$.
\end{itemize}
Condition~(B) arises because $\gamma=\varepsilon\delta$ in Eq.~(30) gives
$\delta\overline{\delta}=\varepsilon^{-1}\overline{\varepsilon}^{-1}$. 
Condition~(A) is tantamount to saying that the ideals $\fr a$ and $\overline{\fr
a}$ are equivalent, in other words that $\fr a$ belongs to a class in $H_1$, the
subgroup of $H\/$ consisting of classes $C\/$ with $\overline{C}=C$.  It is
easily checked that Condition~(B) also depends only on the class of $\fr a$ and
is preserved under multiplication and inversion of classes.  For Condition~(B)
to be applicable at all $\fr a$ must belong to a class in $H_1$.  Consequently
Condition~(B) is equivalent to $\fr a$ belonging to a class in a certain
subgroup $H_2$ of $H_1$.

When Condition~(B) is satisfied the numbers $\gamma=\xi\delta\varepsilon^{-1}$,
where $\xi$ runs through the $N\/$ roots of 1 in $K$, satisfy $|\gamma|=1$.
In this case $\num(\gamma)$ is the ideal ${\fr a}$ and can still be defined
exactly as in Eq.~(\ref{den1}), provided that ``gcd" is interpreted as meaning
``the ideal generated by".  Again ${\fr c}\cap\gamma{\fr c}=\num(\gamma){\fr c}$
and the coincidence index associated with the rotation $\gamma$ is
$\N(\num(\gamma))$ (independent of the ideal $\fr c$).  In the general case,
when the $\omega$'s may be non-principal ideals, a member of the product group
on the right of Eq.~(31) is a pair (root of unity, fractional ideal of the form
${\fr a}\overline{\fr a}^{-1}$) and the argument of Section~3 shows that
elements of $SOC({\fr c})$ correspond precisely to those pairs with the class of
$\fr a$ in $H_2$.  Such pairs form a subgroup of finite index in the full
product group.  One can choose a set of generators for this subgroup in much the
same way as one chooses a basis for a lattice of finite index in a given
lattice, and as in that case there is an infinite number of such bases and no
canonical choice.

Although $SOC({\fr c})$ has independent generators as a group, the set of
coincidence indices in general no longer has independent generators as a
semigroup.

\subsubsection{Coincidence reflections in the general case}

Our second problem was how to calculate indices of coincidence reflections for
a module class in which no module is invariant under complex conjugation. 
Choose, for simplicity, a prime ideal $\fr p$ in the class (which is possible
since every ideal class is known to contain infinitely many prime ideals).  Then
${\fr p}\cap\overline{\fr p}={\fr p}\overline{\fr p}$ has index $\N({\fr p})$
in $\fr p$.  Every coincidence reflection of $\fr p$ has the form
$\rho=\gamma\overline{\,\cdot\,}$ for some $\gamma\in \CC$ with $|\gamma|=1$.
Being a coincidence reflection on $\fr p$, $\rho(\alpha)=\beta$ for some
$\alpha,\beta\in\fr p$.  Hence $\gamma=\beta/\overline{\alpha}\in K$. The index
of $\rho$ is the index of ${\fr p}\cap\gamma\overline{\fr p}$ in $\fr p$ which
is
\be \label{reflec index}
\begin{cases}\N(\num(\gamma))\N({\fr p}),&
   \text{if ${\fr p}\not\hspace*{0.9mm}\mid\num(\gamma)$}, \\
       \N(\num(\gamma))/\N({\fr p}),& 
   \text{if ${\fr p}\mid\num(\gamma)$}.\end{cases}
\ee
We note that there is a reflection of index 1 if and only if the class of 
$\fr p$ is in $H_2$ (when we can choose $\gamma$ to be a generator of the
fractional ideal ${\fr p}/\overline{\fr p}$) and that in that case
(\ref{reflec index}) agrees with our previous way of calculating the index. 
When the class of $\fr p$ is not in $H_2$ the smallest reflection index is
got by taking $(\gamma)={\fr p\,a}/\overline{\fr p}\,\overline{\fr a}$, where
$\fr a$ is the ideal of minimal norm such that the class of $\fr pa$ is in
$H_2$.  Of course, $OC(\overline{\fr p})$ is $OC({\fr p})$ conjugated by
reflection in the $x\/$-axis (corresponding isometries having the same index). 
We note that this is consistent with (\ref{reflec index}): just replace $\fr p$
and $\gamma$ by their complex conjugates.
\begin{theorem}
The group of coincidence rotations of a general $n$-fold
symmetric module is the direct product of its finite rotation symmetry group
$C_N$ and countably many infinite cyclic groups which can be effectively
computed and depend only on $n$.
The index of any coincidence rotation so presented can be calculated explicitly.
Any such module is equivalent to some prime ideal in the cyclotomic field of
$n$th roots of unity, and in this form complex conjugation represents the
coset of coincidence reflections whose indices can be computed from
(\ref{reflec index}) (they depend not only on $n$ but on the individual module).
Such a module need not have exact reflection symmetry.
\end{theorem}

The following table list some statistics for the first few cyclotomic fields
with $h>1$.  We follow Washington$^{28}$ in listing fields with their degree,
$\phi(n)$, as the primary order and $n\/$ as the secondary order.  For each
field we give $n$, $N$, $H$, $H_1$, $H_2$, the smallest rotation index and
the smallest reflection index of the non-principal modules (for the principal
module it is always 1).  In brackets after each index we give the number of
different rotations or reflections with that index.  For all fields on our list
$H_2$ is the trivial subgroup consisting only of the identity element $E\/$ of
$H$.  Also complex conjugation acts on the class group as multiplicative
inversion for all these fields.
$$\begin{array}{ccccccrlrl}
\hline\hline
&&&&&&\multicolumn{2}{c}{\mbox{Min. rotation}}
&\multicolumn{2}{l}{\mbox{Min. reflection index of}}\\
n&N&\mbox{Degree}&H&H_1&H_2&\multicolumn{2}{c}{\mbox{index}}&
\multicolumn{2}{l}{\mbox{non-principal modules}}\\
\hline
23&46&22&C_3&\{E\}&\{E\}&\hspace*{.5cm}599&(22)&\hspace*{1.5cm}47&(11)\\
39&78&24&C_2&C_2&\{E\}&157&(24)&13&(2)\\
52&52&24&C_3&\{E\}&\{E\}&313&(24)&13&(1)\\
56&56&24&C_2&C_2&\{E\}&64&(2)&8&(2)\\
72&72&24&C_3&\{E\}&\{E\}&729&(2)&9&(1)\\
29&58&28&C_2^3&C_2^3&\{E\}&4931&(28)&59&(4)\\
31&62&30&C_9&\{E\}&\{E\}&5953&(30)&
\left\{\begin{array}{r}\phantom{1}32\\125 \end{array} 
\hspace*{-1ex} \right.&
\parbox{3cm}{(1): order 9\\(5): order 3}\\
\hline\hline
\end{array}$$
The two sets of figures in the last entry are due to the fact that
non-principal modules with different orders in the class group of
$\QQ(e^{2\pi i/31})$ have different minimum reflection indices.

\subsubsection{Another example: $N=46$}

To illustrate the results of the previous subsection we treat in detail the case
$n=23$ (with 46-fold symmetry).  For this $n$, the class group $H\/$ of $K\/$ is
$H=\{E,C,C^2\}$, where $C^3=E\/$ and $\overline C=C^2$.  Hence $H_1=\{E\}$ and
therefore $H_2=\{E\}$ too.

The methods of Section~3 show that the complex splitting primes are precisely
those that are quadratic residues mod 23 and for these $\deg(p)=1$ or 11
according to whether $p\equiv1$ mod 23 or not.  The prime ideals ${\fr p}$ of
${\cal O}_K$ that divide a given rational prime $p\/$ are either all principal
or all non-principal (because the Galois group $\Gal(K/\QQ)$ permutes them
transitively) and in the non-principal case fall into complex conjugate pairs of
ideals, one from each of the classes
$C\/$ and $C^2$.  We partition the set of pairs $\Omega$ into the sets
$\Omega_1$, $\Omega_2$ as follows:
\bea
\Omega_1&=&\bigl\{
\{\omega_1,\overline\omega_1\},\{\omega_2,\overline\omega_2\},\ldots
\bigr\}\\
\Omega_2&=&\bigl\{
\{{\fr p}_1,\overline{\fr p}_1\},\{{\fr p}_2,\overline{\fr p}_2\},\ldots
\bigr\}\, ,
\eea
where the $\omega_i$'s are numbers (corresponding to principal ideals) and where
in $\Omega_2$ we have chosen ${\fr p}_i\in C$, $\overline{\fr p}_i\in C^2$
for each $i$.  Finding all numbers of $K\/$ on the unit circle is equivalent to
finding all principal ideals with $K/L$\/-norm equal to ${\cal O}_L$.  (The
numbers $\gamma$ are then the sets of associates of the generators of these
ideals.)  These ideals are precisely those of the form
\be \label{nonprinc generators}
\prod_l\left( \frac{\omega_l}{\overline\omega_l} \right)^{m_l}
\prod_k\left( \frac{{\fr p}_k}{\overline{\fr p}_k} \right)^{n_k}
\ee
with $\sum n_k$ divisible by 3 (since each ${\fr p}_k^{}{\fr p}_k^{-1}$ belongs to the
class $C^2$ of order 3).  This group of ideals has each
$\omega_l/\overline\omega_l$ as an independent generator of the first factor,
and a set of independent generators of the second factor can be chosen as
follows:
\be
({\fr p}_1^{}/\overline{\fr p}_1^{})^3,\quad
\overline{\fr p}_1^{}{\fr p}_2^{}/{\fr p}_1^{}\overline{\fr p}_2^{},\quad
\overline{\fr p}_2^{}{\fr p}_3^{}/{\fr p}_2^{}\overline{\fr p}_3^{},\;\ldots.
\ee

Although this exhibits $SOC({\fr c})$ as having independent generators as a
group, the set of coincidence indices no longer has independent generators as a
semigroup.  Instead of basic coincidence indices one has the prime powers
$$p\quad(p\equiv1\mbox{ (mod 23))\quad and\quad}
p^{11}\quad(p\equiv2,3,4,6,8,9,12,13,16,18\mbox{ (mod 23))}$$
which can be partitioned into two classes $P_1$ and $P_2$ (corresponding to
$\Omega_1$ and $\Omega_2$) according to whether or not the prime ideals dividing
$p\/$ are principal.  
As examples:
$$599,691,829,59^{11},101^{11}\in P_1$$
and
$$47,139,277,461,967,
2^{11},3^{11},13^{11},29^{11},31^{11},41^{11},71^{11},73^{11}\in P_2.$$
These examples were computed using the observation (derived from the last
paragraph of Chapter 1 of \cite{W}) that $p\/$ factorizes into principal primes
if and only if it factorizes into principal primes in $\QQ(\sqrt{-23})$.  A
necessary and sufficient condition for this is the solubility of the 
Diophantine equation $6x^2+xy+y^2=p$.

The general product of these numbers has the form
\be
m=p_1^{a_1}\cdots p_r^{a_r}(p_{r+1}^{11})^{a_{r+1}}\cdots(p_s^{11})^{a_s}\times
\{P_1\mbox{-factors}\},
\ee
where $p_1^{},\ldots,p_s^{11}$ are in $P_2$ with
$p_1^{},\ldots,p_r^{}\equiv1\pmod{23}$ and
$p_{r+1}^{},\ldots,p_s^{}\not\equiv1\pmod{23}$.  Now, for $k=r+1\ldots s$ define
\be
\epsilon_k=\begin{cases}0& \text{if $3\mid a_k$},\\
		  1& \text{if not}.\end{cases}
\ee
Then $m\/$ is a coincidence index if and only if
\be
a_1+\cdots+a_r+\epsilon_{r+1}+\cdots+\epsilon_s\ne1.
\ee
[The reason for this is that in choosing a principal ideal giving index $m\/$ we
can arrange that $\sum n_k$ is divisible by 3 in (\ref{nonprinc generators}) by
changing the sign of some $n_k$'s provided at least two $n_k$'s are not divisible
by 3.  Primes of degree 11 are divisible by only one pair of primes in $K$, but
primes of degree 1 are divisible by 11 such pairs, so for these we can easily
arrange that no $n_k$ is divisible by 3.]

Consequently the first three rotation coincidence indices are 1, 599, 691,
the smallest not composed entirely of primes $\equiv1$ (mod 23) is
$2^{11}47=96256$ and the smallest with no prime factors
$\equiv1$ (mod 23) is $2^{11}3^{11}=362797056$.

The Dirichlet series generating function of $f(m)$ can be found much as before,
except that the contribution from non-principal ideals with $K/L\/$-norm equal
to ${\cal O}_L$ must be omitted.  This can be done using the three characters
of the class group: we form three Dirichlet series (Hecke $L\/$-series), one for
each character, by multiplying each norm in the series by the value of the
character on its ideal. The required generating function is then the average of
these three series.

For the principal character (identically equal to 1) the corresponding
Dirichlet series is exactly as in Eq.~(\ref{Dirichlet series}), namely
\be
\left(1+\frac{1}{23^s}\right)^{-1} \frac{\zeta_K(s)}{\zeta_L(2s)} \; .
\ee
For a non-principal character $\chi$ the Euler factor for a prime $p\/$
occurring in $P_1$ is exactly as in (\ref{Euler long}) and (\ref{Euler short}).
For a prime $p\/$ occurring in $P_2$, however, the Euler factor is
\be
\left(\cdots + \frac{\eta}{p^{2ds}}+\frac{\eta^2}{p^{ds}}+ 1+
\frac{\eta}{p^{ds}}+\frac{\eta^2}{p^{2ds}} + \cdots\right)_{ }^{g/2}=
\left( \frac{1-2p^{-ds}}{1-p^{-ds}}\right)_{ }^{g/2},
\ee
where $\eta^3=1$.  Since this does not depend on which primitive cube root of unity
$\eta$ is, the Dirichlet series formed with the characters $\chi$ and
$\overline{\chi}$ are the same and we have
\bd
\sum_{m=1}^\infty \frac{f(m)}{m^s}=
   \frac{1}{3}\left(1+ \frac{1}{23^s}\right)^{-1}
   \frac{\zeta_K^{}(s)}{\zeta_L^{}(2s)}
   \left\{1+2\prod_{p\in P_2}
   \left( \frac{1-2p^{-s}}{1+p^{-s}} \right)^{11}
   \prod_{q\in P_2}\left( \frac{1-2q^{-s}}{1+q^{-s}} \right)\right\},
\ed
where the first product is over the primes $p\/$ in $P_2$ (which are
$\equiv1$ (mod 23)) and the second is over the 11th powers $q\/$ in $P_2$.

In line with earlier examples we give the first 12 nonzero terms:
{\footnotesize \bd
   1+\frac{22}{599^s}+\frac{22}{691^s}+\frac{22}{829^s}+\frac{22}{1151^s}+
     \frac{110}{2209^s}+\frac{22}{2347^s}+\frac{22}{2393^s}+
     \frac{22}{3037^s}+\frac{22}{3313^s}+\frac{22}{3359^s}+
     \frac{22}{4463^s}+\cdots  
\ed }
Note that this applies to all three modules.

For the principal module the reflection indices are the same as the rotation
indices.  For the two non-principal modules, however, the first three
reflection indices are 47, 139, 277.

\section{Concluding remarks}

Let us summarize our results. We have solved the coincidence problem
for planar patterns with $N$-fold symmetry by number theoretic methods.
The first stage consisted of the analysis of lattices and modules in the
plane where the coincidence indices are integers.

{}For various cases of interest we have given the solution explicitly, in
particular describing the set of possible coincidence indices and the number
of coincidence isometries with given index.  The method is described in
sufficient detail to allow other examples along these lines to be worked out.
This is relatively easy for $N<46$, but the complication increases
astronomically for larger $N$ as foreshadowed even in the example $N=46$,
where the class number is only 3.

The second stage was the explicit
investigation of discrete structures associated with a given
module.  Here, in the non-periodic case, the
calculation of the coincidence ratio requires a non-integral
correction factor. We have demonstrated its calculation
in several examples.

{}Furthermore, the approach via algebraic number fields automatically
yields sets of independent generators for the CSM group and therefore
an explicit description of it. The group structure is interesting
in itself because we deal here with infinite discrete groups that are
countably generated and the structure of such groups is not at
all obvious. 

An obvious next step is to extend the investigation
to 3D examples. This is not only an interesting extension of the
technique, but may have concrete realizations. There are two cases
to consider: first the T-phases, i.e. quasicrystals which have a unique
quasiperiodic plane and are periodic in the third direction. The
CSMs for rotations around the unique axis are the ones treated in this
paper. CSMs around other axes occur only when special relations hold
between the lattice constants in the plane and perpendicular to the
plane, a result familiar from the hexagonal case \cite{grim}. There
are also near-coincidences with small misfits between the two grains,
but it is beyond our scope to deal with these. The
second case is the icosahedral one, the only remaining
non-crystalline symmetry in 3D. 
Here we do not have such a powerful tool as the complex numbers and
the structure of the CSM groups is more complicated, even the
rotation part being non-Abelian in general.
Some results are reported in \cite{BP} and will be described more fully in
\cite{next}.

\section*{Appendix A: other rotation centres}

\addtocounter{equation}{-75}
\renewcommand{\theequation}{A\arabic{equation}}

In the main text, we have analyzed the standard situation of coincidence rotations
around lattice (or module) points. Here, we will briefly comment on rotations
around other centres in the lattice case and on situations with more than one
translation class of points.

\subsubsection*{n=4: the square lattice revisited}

Another obvious rotation problem is that around the centre of a Delaunay cell of
$\ZZ^2$, $(\frac{1}{2},\frac{1}{2})$ say. This point represents the only class of deep
holes of $\ZZ^2$, cf.\ \cite{Conway}, and has the entire point group 
$D_4$ of $\ZZ^2$ as site symmetry. 
It is obvious that the coincidence problem is equivalent to
that of the point set $\Gamma$ defined by
\be \label{shift1}
     \Gamma \, = \, \{ a + i b \; | \;  a,b \in \ZZ \, , \; a+b \mbox{ odd} \} \; ,
\ee
which is obtained from $\ZZ^2 - (\frac{1}{2},\frac{1}{2})$ via rotation through
$\pi/4$ and dilation by $\sqrt{2}$.

Observe that (\ref{shift1}) can be rewritten as
\be \label{shift11}
     \Gamma \, = \, \{ \alpha \in \ZZ[i] \; | \:
              \alpha \equiv \!\!\!\!\! / \;\, 0 \;\; (1\!+\!i) \} 
\ee
which solves the problem: as was shown in Section 2, the coincidence rotations of
$\ZZ^2$ can be factorized, the generators being
$e^{i\varphi}=i$ (rotation through $\pi/2$) or of the form 
$e^{i\varphi} = \omega_p/\overline{\omega}_p$ with $N(\omega_p)=p \equiv 1 \; (4)$,
hence $\omega_p  \equiv \!\!\!\!\!\! / \;\, 0 \; (1+i)$. 
The former still is a symmetry of $\Gamma$ (index 1), and we also get the latter
because both numerator and denominator are in $\Gamma$.
Also, the reflection in the $x$-axis remains a coincidence 
operation of index 1.  Summarizing:
\be \label{shift12}
    OC(\Gamma) \; = \; OC(\ZZ^2) \; ,
\ee
and the coincidence indices are unchanged.

\subsubsection*{n=3: the hexagonal packing}

Consider the Voronoi complex of the triangular lattice --- it is a packing made from
regular hexagons --- and let $H$ be its vertex set. 
Let us consider rotations around the centre of a hexagon
which is a point of maximal site symmetry $D_6$.
If we rotate the complex through $\pi/6$ and dilate by $\sqrt3$, 
then $H$ can be characterized as
\be \label{shift2}
    H \, = \, \{ \alpha \in \ZZ[\varrho] \; | \; 
          \alpha \equiv \!\!\!\!\! / \;\, 0 \;\; (1\!+\!\varrho) \} 
\ee
where $\varrho = \frac{1 + i \sqrt{3}}{2}$.
Since $N(1+\varrho) = 3$ and 3 is not a complex splitting prime in $\QQ(\varrho)$,
we find again all rotations and reflections which we had
already for the triangular lattice:
\be \label{shift21}
    OC(H) \, = \, OC(A_2) \; ,
\ee
and also the indices remain unchanged.

\subsubsection*{n=3: coincidence definition revisited}

Slightly different is the situation if we keep the entire set of lattice points,
but rotate around the centre of a Delaunay cell: the latter is a triangle and
its centre has only $D_3$ site symmetry. We rotate again
through $\pi/6$ and dilate by $\sqrt3$ which gives the point set
\be \label{shift22}
     G \, = \, \{ \alpha \in \ZZ[\varrho] \; | \; 
          \alpha \equiv  1 \;\; (1\!+\!\varrho) \} \; .
\ee
Here, a rotation through $\pi/3$ would change the congruence class of $G$
from $1$ to $-1$, so it is no longer a coincidence rotation. This
reduces the torsion part of $OC$ from $C_6$ to $C_3$ in agreement with the reduced
site symmetry while all other generators remain unchanged.
In particular, the reflection in the $x$-axis leaves $G$ invariant and
the index formula applies for all remaining elements.

One might also consider possible variants of the
coincidence concept here: a rotation through $\pi/3$ alone does not produce
a coincidence for the set $G$, while the same rotation followed by a 
suitable translation can give a coincidence of index 1. The latter might
be more important when the connection to grain boundary growth is
considered.
Indeed, especially in view of applications to nonperiodic discrete
point sets, one might define (with obvious meaning)
\be \label{newdef}
   \inf_{t \in \RR^2} [P \; : \; P \cap (RP+t)]
\ee
to be the coincidence index of an isometry $R$ acting on a point set $P$.
This gets rid of the dependence of the index on the rotation centre
and comes closer to the idea of optimal fitting of grain fragments.

\subsubsection*{n=5: the rhombic Penrose tiling}

A complication here is that the vertex sites of the rhombic Penrose
tiling ${\cal T}$ fall into 4 different translation classes with respect to the 
uniquely defined limit translation module ${\cal M}({\cal T})$, 
compare \cite{BKSZ,Martin}.
We identify ${\cal M}({\cal T})$ with the projection of the 
4D root lattice $A_4^{ }$ into tiling space for definiteness. Then
each point class has its own window of pentagonal shape.
The windows come in two different sizes (related by a factor of 
$\tau = (1+\sqrt{5})/2$) and in pairs related by 
rotation through $\pi$, compare \cite{BKSZ}.
The vertices of the rhombi are not points of the module ${\cal M}({\cal T})$
(which also means that none of them is a ``standard'' rotation centre).  

Let us now consider the coincidence problem of the set of vertex sites with
all translation classes identified. To be explicit, we take the rhombic version
of the cartwheel pattern where the rotation centre is not a rhombus vertex but
coincides with the centre of a regular decagon filled with rhombi.
This point is a representative of the fifth translation class,
so far absent.  The cartwheel tiling has $D_{10}$ symmetry in the sense that
any $D_{10}$-operation either maps the tiling upon itself (thus, in particular,
the set of vertex sites) or produces at most a mismatch of density zero (along
the well-known worms).  All these operations thus have coincidence ratio~1. The
corresponding rotation in window space maps windows to windows, because they
appear in $D_{10}$-orbits around the origin.  More than this, it maps
translation classes of windows to translation classes of windows.

For other coincidence isometries, we first observe that the the integral span of
all vertex points is again a planar module of rank 4, in our explicit case the
projection of the weight lattice $A_4^*$, the dual of $A_4^{ }$, into tiling space.
This module is equivalent to ${\cal M}({\cal T})$
and possesses therefore the same coincidence isometries, namely those described
in Section 3. Consequently, we find all these also as coincidence isometries of
the rhombic cartwheel tiling. The coincidence ratio must now be corrected in
a similar way
to that of the Ammann--Beenker tiling in Section 4, but the window system
requires a slightly more complicated calculation which we will not present here.

Even more complicated would be the coincidence analysis for rotations around vertex
points, in particular with various point classes distinguished. The methods needed
are in principle those described for $n=3$ above, but details will not be given
here.

\subsubsection*{n=12: a square-triangle tiling}

Quasiperiodic square-triangle tilings are attractive for a number of reasons.
We mention them because they can have 12-fold symmetry in the sense of
mismatches of at most density zero under $D_{12}$-operations, see \cite{BKS}
for an example. There, all vertex points are in one translation class,
so no problem occurs and we find all coincidence isometries of Section 3.
But for the correction factor due to window overlaps
one encounters a new type of complication: the window is fractally shaped
and consequently we see no way of calculating this factor.  It is left as an
exercise for fractal readers.

\section*{Appendix B: proofs}

\addtocounter{equation}{-7}
\renewcommand{\theequation}{B\arabic{equation}}

Here we give the promised references and proofs of Facts 1--3 in
Section~\ref{CN1}.

{\bf Fact 1}\quad This is proved in \cite{CF} (Lemma 4) or
\cite{W} (Theorem 2.13) for example, but we sketch the proof
here as it leads on naturally to the proof of Fact 2, which is less
commonly found in the literature.

Let ${\fr P}$ be a prime factor of $p\/$ in $K$.
Since $p \nmid n\/$ the $n\/$th roots of 1 
in $K\/$ are distinct mod ${\fr P}$.  (The
most straightforward way to see this is from the identity
\be
n=\prod_{k=1}^{n-1}(1-\xi^k),
\ee
got by putting $x=1$ in $(x^n-1)/(x-1)$, and noting that every difference of
roots of unity is an associate of $1-\xi^k$ for some $k$.)   
The residue class field
$\FF^{}_{\frs P}=\ZZ[\xi]/{\fr P}$ is a finite field generated over
$\FF^{}_p$ by
the residue class $\xi^*$ of $\xi$, and since distinct roots of unity are
distinct mod $\fr P$, the order of $\xi^*$ in $\FF^{}_{\frs P}$ is $n$.  Every
finite extension of $\FF^{}_p$ is
normal with cyclic Galois group generated by the {\em Frobenius automorphism}
$x\mapsto x^p$, whose order is the degree $d\/$ of the extension.
Consequently the degree $[\FF^{}_{\frs P}:\FF^{}_p]$ is the smallest $d\/$ with
$\xi^{*p^d}=\xi^*$; that is, the smallest $d\/$ with $n\mid(p^d-1)$.  This
establishes Fact 1, since $[\FF^{}_{\frs P}:\FF^{}_p]$ is the degree of the
minimal polynomial satisfied by $\xi$ mod $p$ and hence is $\deg^{}_K(p)$.

{\bf Fact 2}\quad Analogously to the above proof, $\deg^{}_L(p)$ is the degree
$d'$ of the residue class field extension $\FF^{}_{{\frs p}}/\FF^{}_p$, where 
${\fr p}$ is the
prime of $L\/$ divisible by ${\fr P}$ and $\FF^{}_{\frs p}$ is the residue class field
$\ZZ[\xi+\xi_{}^{-1}]/{\fr p}$.  
Clearly $[\FF^{}_{\frs P}:\FF^{}_{\frs p}]\le[K:L]=2$, so the
$d\/$ of Fact 1 is either $d'$ or $2d'$.  If $d\/$ is {\em odd} then $d'=d$, the
order of $p\/$ mod $n$, and no power of $p\/$ is congruent to $-1$ mod $n$.

To treat the case of {\em even} $d\/$ we first note that if
$\xi^{*k_1}+\xi^{*-k_1}=\xi^{*k_2}+\xi^{*-k_2}$ in $\FF^{}_{\frs P}$
(where $0\le k_1,k_2<n$) then
either $\xi^{*k_1}=\xi^{*k_2}$ or $\xi^{*k_1}=\xi^{*-k_2}$.  This is because
$\xi^{*k_j}$, $\xi^{*-k_j}$ are the two roots in $x\/$ of
\be
x^2-(\xi^{*k_j}+\xi^{*-k_j})x+1=0,\quad(j=1,2)
\ee
and when the equations are the same the roots must match in some order.  Now
$d'=d/2$ if and only if $\FF^{}_{\frs p}$ is the unique subfield of index 2 in
$\FF^{}_{\frs P}$, this being the fixed field of the element $x\mapsto x^{p^{d/2}}$
of order 2 in the Galois group of $\FF^{}_{\frs P}/\FF^{}_p$.  So $d'=d/2$ if and only
if
\be
\xi^*+\xi^{*-1}=(\xi^*+\xi^{*-1})^{p^{d/2}}=\xi^{*p^{d/2}}+\xi^{*-p^{d/2}},
\ee
which requires $\xi^{*-1}=\xi^{*p^{d/2}}$ (equivalent to $n\mid p^{d/2}+1$),
since $\xi^*\ne\xi^{*p^{d/2}}$.  The exponent $d/2$ here is plainly minimal,
since $n\mid p^a+1\Rightarrow n\mid p^{2a}-1$.

{\bf Fact 3}\quad Part (a) is a result of the fact that when $n=p^r$ then $p\/$
is totally ramified in $K\/$ (that is, $p\/$ is the $\phi(p^r)$th power of a
degree~1 prime of $K\/$), see for example \cite{CF}, Lemma 3.  As a consequence
$\deg^{}_K(p)=1$ and $p\/$ is not a complex splitting prime because it has only
one prime factor in $K$.

\begin{figure}
\centerline{\includegraphics{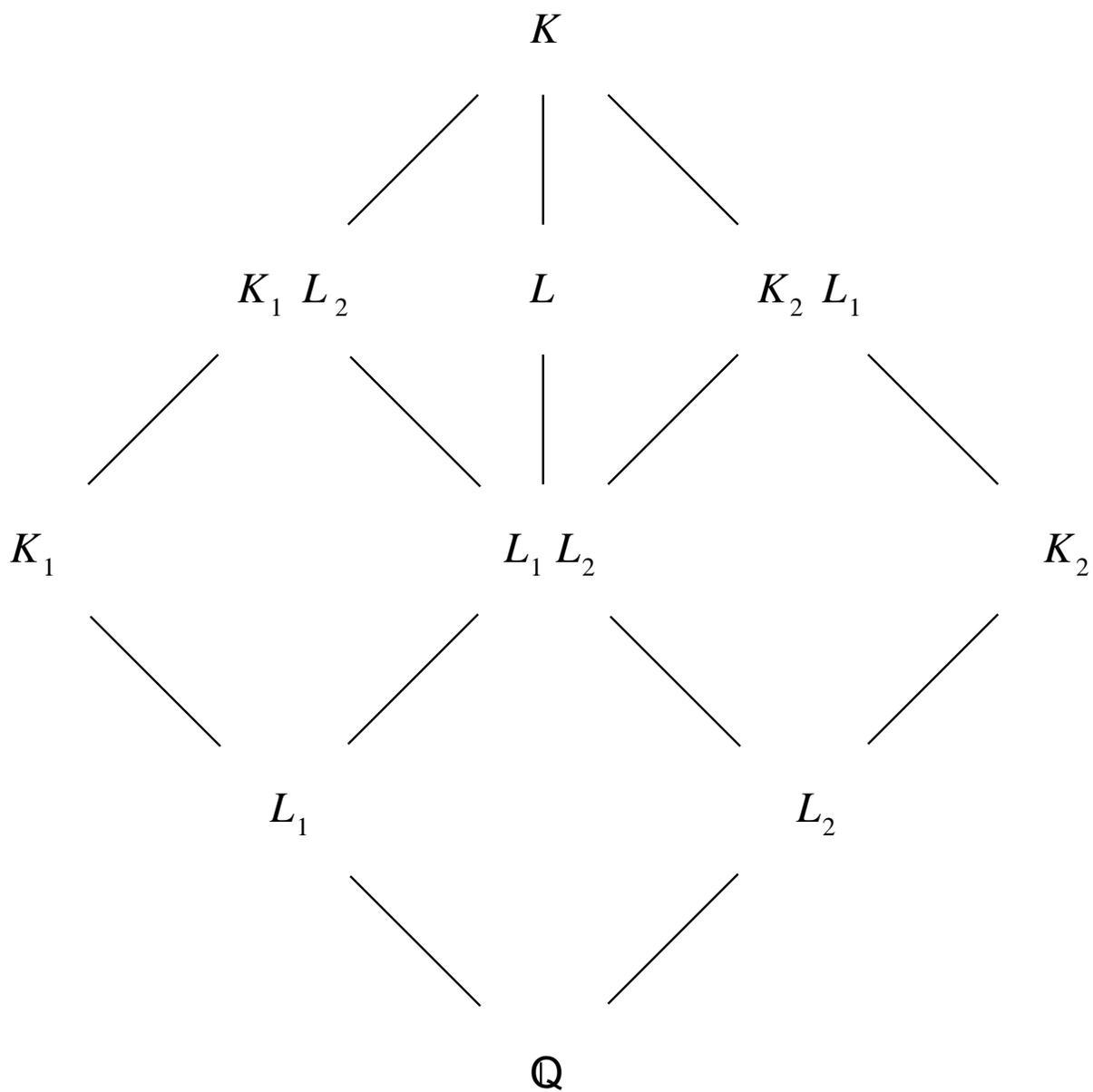}}
\caption{Hasse diagram.}
\label{fig:hasse}
\end{figure}

For part (b) we refer to the Hasse diagram of field inclusions in Figure 4.
Here $K_1$ and $K_2$ are the cyclotomic fields of $n_1$th and
$p^r$th roots of unity and $L_1$ and $L_2$ their maximal real subfields.
Then $K=K_1K_2$, the compositum of $K_1$ and $K_2$, and, since
$p\not\hspace*{0.9mm}\mid n_1$, $K_1\cap K_2=\QQ$ (see \cite{Cohn} Thm.~9.52 or \cite{W}
Prop.~2.4).  Let ${\fr p}={\fr p}^{}_K$ be a prime of
$K$ dividing $p$.  For an arbitrary subfield $F\/$ of $K\/$ we denote by
${\fr p}^{}_F$ the prime ideal of $F\/$ that is divisible by ${\fr p}$.

Because $p\/$ is unramified in $L_1$ and $K_1$ but totally ramified in $K_2$ it
follows that ${\fr p}^{}_{L_1}$ and ${\fr p}^{}_{K_1}$ are totally ramified in
$K_2L_1$ and $K\/$ and, in particular,
\be
\deg^{}_{K_2L_1/L_1}({\fr p}^{}_{K_2L_1})=\deg^{}_{K/K_1}({\fr p}^{}_{K})=1.
\ee
Consequently
\be \label{equaldeg1}
\deg^{}_{K/K_2L_1}({\fr p})=\deg^{}_{K_1/L_1}({\fr p}^{}_{K_1}).
\ee
Now look at the fields $L_1L_2$, $K_2L_1$, $L\/$ and $K\/$.  Since
${\fr p}^{}_{L_2}$ ramifies in $K_2$ but $p\/$ is unramified in $L_1$,
${\fr p}^{}_{L_1L_2}$ ramifies in $K_2L_1$ and hence in $K$.  By Prop.~2.15(b) of
\cite{W}, ${\fr p}^{}_L$ is unramified in $K$, and hence ${\fr p}^{}_{L_1L_2}$
ramifies in $L$.  We now have
\be
\deg^{}_{L/L_1L_2}({\fr p}^{}_L)=\deg^{}_{K_2L_1/L_1L_2}({\fr p}^{}_{K_2L_1})=1
\ee
whence
\be \label{equaldeg2}
\deg^{}_{K/L}({\fr p})=\deg^{}_{K/K_2L_1}({\fr p}).
\ee
Since ${\fr p}^{}_L$ is unramified in $K\/$ and ${\fr p}^{}_{L_1}$ is unramified in $K_1$,
Eqs.~(\ref{equaldeg1}) and (\ref{equaldeg2}) imply that ${\fr p}^{}_L$ factors into two
primes of $K\/$ if and only if ${\fr p}^{}_{L_1}$ factors into two primes of $K_1$.
Finally, $\deg^{}_K(p) = \deg^{}_{K_1}(p)$ is an immediate consequence
of the fact that ${\fr p}_{K_1}$ is totally ramified in $K$.

\vspace{5mm}
\parindent 0pt
{\bf Acknowledgements}
\vspace{2mm}
\parindent 15pt

It is a pleasure to thank Reinhard L\"uck and David Warrington for getting us
interested in this problem and for many stimulating discussions in the early
stages of this work. 
We are grateful to Nicolas Rivier for many valuable hints and
the second referee for various corrections.
P.~P.\ is grateful to the Universities of Stuttgart and T\"ubingen for
hospitality while this work was carried out and to Pieter Moree for helpful
suggestions.

\end{document}